\documentclass[11pt,fleqn]{article}

\usepackage{amsfonts,amssymb,amsthm,latexsym,mathrsfs}

\newcommand{\dd}{{\rm d}}
\newcommand{\ee}{{\rm e}}
\newcommand{\oo}{{\rm o}}

\newcommand{\OO}{{\rm O}}

\newcommand{\AAAA}{{\bf A}}
\newcommand{\BB}{{\bf B}}
\newcommand{\DD}{{\bf D}}
\newcommand{\II}{{\bf I}}
\newcommand{\UU}{{\bf U}}
\newcommand{\VV}{{\bf V}}
\newcommand{\WW}{{\bf W}}
\newcommand{\XX}{{\bf X}}

\newcommand{\expect}[1]{{\mathbb E}\left\{#1\right\}}
\newcommand{\pr}[1]{{\mathbb P}\left\{#1\right\}}

\newtheorem{theorem}{Theorem}
\newtheorem{lemma}[theorem]{Lemma}
\newtheorem{proposition}[theorem]{Proposition}
\newtheorem{example}[theorem]{Example}

\begin{document}

\title{\mbox{Heavy Loads and Heavy Tails}}

\author{Sem Borst \\
Department of Mathematics \& Computer Science \\
Eindhoven University of Technology \\
P.O. Box 513, 5600 MB Eindhoven, The Netherlands}

\maketitle

\begin{abstract}
The present paper is concerned with the stationary workload of queues with heavy-tailed (regularly varying) characteristics.
We adopt a transform perspective to illuminate a close connection between the tail asymptotics and heavy-traffic limit in infinite-variance scenarios.
This serves as a tribute to some of the pioneering results of J.W. Cohen in this domain.
We specifically demonstrate that reduced-load equivalence properties established for the tail asymptotics of the workload naturally extend to the heavy-traffic limit.
%In contrast, while the tail asymptotics are governed by the service requirement distribution and not affected by the interarrival time distribution, it turns out that the interarrival times and service requirements play a symmetric role for the heavy-traffic limit.

\end{abstract}

\section{Introduction}

Queueing theory has primarily been concerned with light-tailed input processes, meaning that the interarrival and service requirement distributions exhibit exponentially decaying tails, and thus have finite moments.
In particular, heavy-traffic analysis of queues, which was initiated by Kingman in the early sixties \cite{Kingman61,Kingman62} and has burgeoned into a huge branch of the queueing literature since,
%pertaining to the behavior in a critical-load regime,
see for instance \cite{Whitt02,Williams16}, typically focuses on the finite-variance case.

Nevertheless, the last few decades have witnessed an emerging interest in queueing models with heavy-tailed input processes where the involved tail distributions fall off slower than any exponential distribution and may have infinite variance.
One of the earliest results in this domain dates back nearly half a century to the seminal paper by Cohen~\cite{Cohen73} who derived the tail asymptotics of the stationary waiting time in a single-server queue with a so-called regularly varying service requirement distribution, see also Borovkov~\cite{Borovkov76}.
Pakes~\cite{Pakes75} and Veraverbeke~\cite{Veraverbeke77} extended these results to the supremum of random walks with increments with a subexponential distribution, see also Embrechts \& Veraverbeke~\cite{EV82}.
Many of these threads were pursued in the context of the ruin probability in closely related insurance risk models, motivated by the observation that claim sizes commonly exhibit power law characteristics.

The interest in queueing models with heavy-tailed input processes received a strong impetus in the late nineties when measurement studies revealed that traffic processes in packet-switched data communication networks exhibit long-range dependence and self-similarity \cite{BSTW95,LTWW94,PF95,WTLW95}.
These findings provided a stark contrast with the traditional Markovian assumptions for circuit-switched voice systems, and may be attributed to heavy-tailed features in the underlying activity patterns, e.g. session lengths or file sizes with infinite empirical variance \cite{CB96,WTSW97}.
The specific features of data communication networks further gave rise to a wide spectrum of model variants beyond the single-server queue, such as fluid queues fed by On-Off sources with heavy-tailed activity periods, see for instance \cite{AMN99,Boxma96,Boxma97,BD98b,JL99,MB00}.

In the latter context a crucial dichotomy was observed between on the one hand scenarios where just a single extremely large service request is enough to cause a large queue length, waiting time or workload, and on the other hand scenarios where several large service requests must overlap in order for such a rare event to occur, rendering the derivation and form of the tail asymptotics significantly more involved \cite{Zwart01,ZBM04}.
The latter type of scenario also arises in multi-server queueing systems where the required number of large service requirements specifically depends on the slack capacity \cite{FK06,FK12,Whitt00}.
This situation has motivated a strong interest in the broader topic of large-deviations behavior when several `big jumps' are necessary to cause a rare event of interest.
A detailed discussion of the recent developments in this direction as well as extensions to network settings \cite{FM14,LM08}, multi-dimensional processes and semi-exponential (e.g. Weibull) distributions \cite{BBRZ19,BBRZ20,CBRZ19} is beyond the scope of the present paper, and for background information and further references we refer to the excellent books \cite{FKZ13,NWZ22}.

Besides the tail properties, the heavy-traffic behavior of the single-server queue with heavy-tailed input processes has attracted strong interest as well.
Some of the earliest results were established in several deep and only partly published studies by Cohen \cite{Cohen97a,Cohen97b,Cohen98a} using transform methods, and later generalized by Boxma \& Cohen \cite{BC99a,BC99b}.
These results showed that for a certain class of heavy-tailed service requirement distributions with infinite variance, which includes the regularly varying distributions mentioned earlier, a suitably scaled version of the stationary waiting time has a Mittag-Leffler distribution in the limit, provided that the tail of the interarrival time distribution is `less heavy'.
Interestingly, the scaling involves a so-called contraction coefficient that differs from the classical finite-variance case and depends on the degree of `heaviness' of the service requirement distribution.
Related process-level limits were established by Resnick \& Samorodnitsky~\cite{RS00}, and pointed to strong connections with $\alpha$-stable processes.
For the dual case where the interarrival time distribution belongs to the above class of heavy-tailed distributions and the service requirement distribution is `less heavy', the results in \cite{BC99a,BC99b} showed that the scaled version of the stationary waiting time has a negative exponential distribution in the limit, just like in the traditional finite-variance case, except that the scaling involves a similar contraction coefficient that then depends on the `heaviness' of the interarrival time distribution.
The case where the tails of the service requirement and interarrival time distributions are `equally heavy' was also covered in \cite{BC99b,Cohen97c}.

The present paper serves as a testimonial to some of the pioneering contributions of J.W. Cohen in the above-described domain, in particular references \cite{Cohen73,Cohen97c,Cohen98a}.
The primary goal is not to establish entirely novel results, but to present an overarching view on some of the known results for the workload distribution in infinite-variance scenarios.
Specifically, we use a transform perspective to highlight a close connection between the tail asymptotics on the one hand and the heavy-traffic limit on the other hand.
It will further be demonstrated that so-called reduced-load equivalence properties for the tail asymptotics of the workload have a natural counterpart in terms of the heavy-traffic limit.
We will capitalize on available transform results for several specific models, including an M/G/1 queue with alternating service speed and an M/G/2 queue with heterogeneous input as illustrative examples.
These model instances all pertain to the scenario where a single large service request is sufficient to cause a large workload, and do not cover any network settings.
%In contrast, while the possibly heavy-tailed nature of the interarrival time distribution does not matter for the workload asymptotics, it does play a similar role for the heavy-traffic limit as a heavy-tailed service requirement distribution.

The remainder of the paper is organized as follows.
In Section~\ref{modelsetup} we present the model set-up and some preliminaries.
Section~\ref{connection} establishes a connection between the tail asymptotics and heavy-traffic limit for queueing systems with regularly varying characteristics in infinite-variance scenarios, and discusses a few basic illustrative examples.
We use this perspective to explore the relationship between these two limiting regimes in the context of queueing systems with time-varying service speed and server heterogeneity in Sections~\ref{alternating} and~\ref{heterogeneous}, respectively.

\section{Model set-up and preliminaries}
\label{modelsetup}

Let $\VV_\lambda$ be a random variable having the distribution of the stationary workload in some queueing system with arrival rate~$\lambda$.
We assume that the stability condition is $\lambda < \lambda^*$, with $\lambda^* \in (0, \infty)$ representing the `critical' arrival rate.

The \textit{tail asymptotics} of the workload concern the behavior of $\pr{\VV_\lambda > x}$ as $x$ grows large, for a given arrival rate~$\lambda$.
In the context of communication networks this is also commonly referred to as \textit{large-buffer asymptotics}.
For brevity, we will write $\pr{\VV_\lambda > x} \sim f(x)$ as $x \to \infty$ for some function $f(\cdot)$ to denote that $\pr{\VV_\lambda > x} = (1 + \oo(1)) f(x)$ as $x \to \infty$, i.e., $\lim_{x \to \infty} \pr{\VV_\lambda > x} / f(x) = 1$.
We will say that $\VV_\lambda$ has the \textit{same tail asymptotics} as the random variable $\hat{\VV}_{\hat\lambda}$ if $\pr{\VV_\lambda > x} \sim \pr{\hat{\VV}_{\hat\lambda} > x} \sim f(x)$ as $x \to \infty$ for some function $f(\cdot)$.

The \textit{heavy-traffic limit} of the workload relates to the distribution of a properly scaled version of $\VV_\lambda$ as the arrival rate~$\lambda$ tends to the critical value~$\lambda^*$.
For conciseness, we will write $S_\lambda \VV_\lambda \to \WW$ as $\lambda \uparrow \lambda^*$ for some random variable $\WW$, with $S_\lambda$ a suitable scaling factor, to indicate that $S_\lambda \VV_\lambda$ converges to~$\WW$ in distribution as $\lambda \uparrow \lambda^*$, i.e., $\lim_{\lambda \uparrow\lambda^*} \pr{S_\lambda \VV_\lambda \leq x} = \pr{\WW \leq x}$ in every point of continuity of $\pr{\WW \leq x}$.
In case $S_\lambda$ is strictly decreasing as function of~$\lambda$ with $S_\lambda \downarrow 0$ as $\lambda \uparrow \lambda^*$, this may be equivalently expressed as $\lim_{\epsilon \downarrow 0} \pr{\epsilon \VV_{\lambda(\epsilon)} \leq x} = \pr{\WW \leq x}$, with $\lambda(\epsilon) = S_\lambda^{- 1}(\epsilon)$, where we suppress the implicit dependence of the function $\lambda(\epsilon)$ on $S_\lambda$ for ease of notation.
We will say that $\VV_\lambda$ has the \textit{same heavy-traffic limit} for $\lambda \uparrow \lambda^*$ as the random variable $\hat{\VV}_{\hat\lambda}$ for $\hat\lambda \uparrow \hat\lambda^*$, up to a relative capacity slack factor $\zeta \in (0, \infty)$, if $\lim_{\epsilon \downarrow 0} \pr{\epsilon \VV_{\lambda(\epsilon)} \leq x} = \lim_{\epsilon \downarrow 0} \pr{\epsilon \hat{\VV}_{\hat\lambda(\epsilon)} \leq x} = \pr{\WW \leq x}$ in every point of continuity for some proper random variable~$\WW$, with $\lambda(\epsilon) \uparrow \lambda^*$ as $\epsilon \downarrow 0$ and
\[
\hat\lambda(\epsilon) =
\hat\lambda^* \left[1 - \zeta \left(1 - \frac{\lambda(\epsilon)}{\lambda^*}\right)\right],
\]
so that
\[
\frac{\hat\lambda^* - \hat\lambda(\epsilon)}{\hat\lambda^*} =
\zeta \frac{\lambda^* - \lambda(\epsilon)}{\lambda^*}.
\]

We will mainly be interested in scenarios where some random variables (e.g.~service requirements) have regularly varying distributions.
A non-negative random variable $\XX$ is called regularly varying of index~$- \alpha$, $\alpha \geq 0$, if %$\pr{\XX > y x} \sim y^{- \alpha} \pr{\XX > x}$ as $x \to \infty$ for all $y > 0$, or equivalently,
$\pr{\XX > x} \sim L(x) x^{- \alpha}$ as $x \to \infty$, with $L(\cdot)$ a slowly varying function, i.e., $L(y x) \sim L(x)$ as $x \to \infty$ for all $y > 0$.

The next lemma presents a useful relationship between the tail asymptotics of $\pr{\XX > x}$ as $x \to \infty$ and the behavior of its Laplace-Stieltjes Transform (LST) $\expect{\ee^{- \omega \XX}}$ as $\omega \downarrow 0$.
It has been formulated as Lemma 2.2 in~\cite{BD98b} as an extension of Theorem 8.1.6 in~\cite{BGT87}.

\begin{lemma}

\label{tauber}

Let $\XX$ be a non-negative random variable, $L(\cdot)$ a slowly varying function, $\alpha \in (n, n + 1)$ with $n \in \mathbb{N}$ and $C \geq 0$.
Then the following two statements are equivalent:

(i) $\pr{\XX > x} = (C + \oo(1)) L(x) x^{- \alpha}$ as $x \to \infty$;

(ii) $\expect{\XX^n} < \infty$ and
\[
\expect{\ee^{- \omega \XX}} - \sum_{m = 0}^{n} \frac{(- 1)^m \expect{\XX^m} \omega^m}{m!} =
- \Gamma (1 - \alpha) (C + \oo(1)) L(1 / \omega) \omega^\alpha \mbox{ as } \omega \downarrow 0.
\]

\end{lemma}

Denote by $\XX^r$ the \textit{excess} or \textit{residual-lifetime} random variable associated with~$\XX$, i.e., $\pr{\XX^r > x} = \frac{1}{\expect{\XX}} \int_{y = x}^{\infty} \pr{\XX > y} \dd y$ and LST $\expect{\ee^{- \omega \XX^r}} = \frac{1 - \expect{\ee^{- \omega \XX}}}{\expect{\XX} \omega}$.
If $\XX$ is regularly varying of index~$- \nu$, $\nu > 1$, i.e., $\pr{\XX > x} = (C + \oo(1)) L(x) x^{- \nu}$ as $x \to \infty$, then Karamata's theorem implies that $\pr{\XX^r > x} = \frac{C + \oo(1)}{(\nu - 1) \expect{\XX}} L(x) x^{1 - \nu}$ as $x \to \infty$.
In case $\nu \in (1, 2)$ we obtain, applying Lemma~\ref{tauber} with $\alpha = \nu$ and $n = 1$, $\expect{\ee^{- \omega \XX}} - 1 + \expect{\XX} \omega = - \Gamma(1 - \nu) (C + \oo(1)) L(1 / \omega) \omega^\nu$ as $\omega \downarrow 0$, and thus $\expect{\ee^{- \omega \XX^r}} - 1 = \frac{1 - \expect{\ee^{- \omega \XX}}}{\expect{\XX} \omega} - 1 = \frac{\Gamma(1 - \nu) (C + \oo(1))}{\expect{\XX}}
L(1 /\omega) \omega^{\nu - 1}$ as $\omega \downarrow 0$.
Alternatively, this may be deduced from the expression for $\pr{\XX^r > x}$, invoking Lemma~\ref{tauber} with $\alpha = \nu - 1$ and $n = 0$, and noting that $\Gamma(2 - \nu) = (1 - \nu) \Gamma(1 - \nu)$.

In the remainder of the paper we will omit the slowly-varying function to eschew the corresponding notational burden and technical intricacies, thus essentially restricting to Pareto-type distributions.
In this regard it is fitting though to draw attention to the tour de force demonstrated by Cohen in many of his studies in explicitly dealing with regularly varying distributions in full generality, in particular avoiding the need to work with implicitly defined De Bruijn conjugates~\cite{BGT87}.

\vspace*{.2in}

\textbf{Tail asymptotics for single-server queue}

For later reference, let $\VV_{\hat\lambda}^c$ be the stationary workload in a single-server queue with service speed~$c$, a Poisson arrival process of rate~$\hat\lambda$ and service requirements that are i.i.d.~copies of a non-negative random variable~$\hat{\BB}$ with mean $\hat\beta < \infty$ and LST $\hat\beta(\omega)$.
When we consider the amount of work $\UU_{\hat\lambda}^c = \VV_{\hat\lambda}^c / c$ measured in time units, this system can be viewed as an M/G/1 queue where the service requirements are scaled by a factor~$1 / c$.

If $\hat{\BB}$ is regularly varying of index~$- \nu$, $\nu > 1$, i.e.,
\begin{equation}
\pr{\hat{\BB} > x} \sim - \frac{C_{\hat{B}}}{\Gamma(1 - \nu)} x^{- \nu}
\hspace*{.4in} \mbox{ as } x \to \infty,
\label{sere1}
\end{equation}
then it follows from~\cite{Borovkov76,Cohen73} that for any $\hat\lambda < \hat\lambda^*$,
\[
\pr{\UU_{\hat\lambda}^c > y} \sim
\frac{\hat\lambda \hat\beta / c}{1 - \hat\lambda \hat\beta / c} \pr{\hat{\BB}^r / c > y} = \frac{\hat\lambda \hat\beta}{c - \hat\lambda \hat\beta} \pr{\hat{\BB}^r > c y}
\hspace*{.4in} \mbox{ as } y \to \infty,
\]
and hence
\begin{equation}
\pr{\VV_{\hat\lambda}^c > x} =
\pr{\UU_{\hat\lambda}^c > x / c} \sim \frac{\hat\lambda \hat\beta}{c - \hat\lambda \hat\beta} \pr{\hat{\BB}^r > x}
\hspace*{.4in} \mbox{ as } x \to \infty,
\label{tailvhlcx}
\end{equation}
with
\begin{equation}
\pr{\BB^r > x} \sim \frac{C_{\hat{B}}}{(1 - \nu) \Gamma(1 - \nu) \hat\beta} x^{1 - \nu} = \frac{C_{\hat{B}}}{\Gamma(2 - \nu) \hat\beta} x^{1 - \nu}
\hspace*{.4in} \mbox{ as } x \to \infty.
\label{tailbrx}
\end{equation}

\iffalse

It is worth observing that $\expect{\ee^{- \omega \VV_{\hat\lambda}^c}}$ in~(\ref{lstvhlc}) may be written in the form
\[
\expect{\ee^{- \omega \VV_{\hat\lambda}^c}} =
\frac{\hat\lambda^* - \hat\lambda}{\hat\lambda^* - \hat\lambda \hat\beta^r(\omega)} =
\frac{1 - \frac{\hat\lambda \hat\beta}{c}}{1 - \frac{\hat\lambda \hat\beta}{c} \hat\beta^r(\omega)} =
\left(1 - \frac{\hat\lambda \hat\beta}{c}\right) \sum_{n = 0}^{\infty} \left(\frac{\hat\lambda \hat\beta}{c} \hat\beta^r(\omega)\right)^n,
\]
which implies
\[
\pr{\VV_{\hat\lambda}^c > x} = \left(1 - \frac{\hat\lambda \hat\beta}{c}\right) \sum_{n = 0}^{\infty} \left(\frac{\hat\lambda \hat\beta}{c}\right)^n \pr{\XX_1 + \dots + \XX_n > x},
\]
where $\XX_1, \XX_2, \dots$ are i.i.d.\ copies of $\hat{\BB}^r$.
Using the properties of regularly varying (in fact subexponential) distributions and dominated convergence, it then follows that
\[
\pr{\VV_{\hat\lambda}^c > x} \sim \left(1 - \frac{\hat\lambda \hat\beta}{c}\right) \sum_{n = 0}^{\infty} \left(\frac{\hat\lambda \hat\beta}{c}\right)^n n \pr{\hat{\BB}^r > x} = \frac{\hat\lambda \hat\beta}{c - \hat\lambda \hat\beta} \pr{\hat{\BB}^r > x}
\hspace*{.4in} \mbox{ as } x \to \infty,
\]
as stated in~(\ref{tailvhlcx}).

\fi

\vspace*{.2in}

\textbf{Tail asymptotics for fluid queue with On-Off-source}

For later use, we also introduce a related fluid queue with drain rate~$d$ fed by a single On-Off source,
which generates input at some rate $r > d$ when On.
The On-periods are i.i.d.\ copies of a non-negative random variable~$\AAAA$ with mean $\expect{\AAAA} < \infty$.
The Off-periods are i.i.d.\ copies of a non-negative random variable~$\UU$ with mean $\expect{\UU} < \infty$.
Define $p_{On} = \frac{\expect{\AAAA}}{\expect{\AAAA} + \expect{\UU}}$ as the fraction of time that the source is On.
Denote by $\rho = p_{On} r$ the time-average input rate of the source, and assume $\rho < d$ for stability.
Let $\WW^d$ be the stationary workload.

If $\pr{\AAAA > x}$ is regularly varying of index~$- \nu$, $\nu > 1$, i.e.,
\begin{equation}
\pr{\AAAA > x} \sim - \frac{C_A}{\Gamma(1 - \nu)} x^{- \nu}
\hspace*{.4in} \mbox{ as } x \to \infty,
\label{onpe1}
\end{equation}
then it follows from~\cite{JL99} that
\begin{equation}
\hspace*{-.4in}
\pr{\WW^d > x} \sim (1 - p_{On}) \frac{\rho}{d - \rho} \pr{\AAAA^r > \frac{x}{r - d}} =
p_{On} \frac{r - \rho}{d - \rho} \pr{\AAAA^r > \frac{x}{r - d}}
\hspace*{.4in} \mbox{ as } x \to \infty,
\label{tailwcx}
\end{equation}
with
\begin{equation}
\pr{\AAAA^r > x} \sim \frac{C_A}{(1 - \nu) \Gamma(1 - \nu) \expect{\AAAA}} x^{1 - \nu} = \frac{C_A}{\Gamma(2 - \nu) \expect{\AAAA}} x^{1 - \nu}
\hspace*{.4in} \mbox{ as } x \to \infty.
\label{tailarx}
\end{equation}

Comparison of~(\ref{tailwcx},\ref{tailarx}) with (\ref{tailvhlcx},\ref{tailbrx}) shows that, under the assumptions~(\ref{sere1}) and~(\ref{onpe1}), $\VV_{\hat\lambda}^c$ and $\WW^d$ have the same tail asymptotics when
\[
(1 - p_{On}) \frac{\rho}{d - \rho} \frac{C_A}{\expect{\AAAA}} (r - d)^{\nu - 1} =
p_{On} \frac{r - \rho}{d - \rho} \frac{C_A}{\expect{\AAAA}} (r - d)^{\nu - 1} =
%\frac{r - \rho}{d - \rho} \frac{C_A}{\expect{\AAAA} + \expect{\UU}} (r - d)^{\nu - 1} =
\frac{\hat\lambda \hat\beta}{c - \hat\lambda \hat\beta} \frac{C_{\hat{B}}}{\hat\beta} =
\frac{\hat\lambda C_{\hat{B}}}{c - \hat\lambda \hat\beta}.
\]
In case $\hat{\BB} = (r - d) \AAAA$, so that $C_{\hat{B}} = C_A (r - d)^\nu$ and $\hat\beta = (r - d) \expect{\AAAA}$, the above equality is satisfied when
\[
(1 - p_{On}) \frac{\rho}{d - \rho} = \frac{\hat\lambda \hat\beta}{c - \hat\lambda \hat\beta},
\]
or equivalently,
\[
c = \left(1 + \frac{1}{1 - p_{On}} \frac{d - \rho}{\rho}\right) \hat\lambda \hat\beta =
\frac{\rho (1 - p_{On}) + d - \rho}{\rho (1 - p_{On})} \hat\lambda \hat\beta =
\frac{d - p_{On} \rho}{\rho - p_{On} \rho} \hat\lambda \hat\beta.
\]

\section{Tail asymptotics and heavy-traffic limit}
\label{connection}

In this section we focus on the tail asymptotics and heavy-traffic limit of the stationary workload $\VV_\lambda$ in queueing systems with regularly varying characteristics as discussed in the previous section.
In order to describe the connection between these two limiting regimes in infinite-variance scenarios, it will be convenient to write the LST of $\VV_\lambda$ for various specific model instances in the form
\begin{equation}
\expect{\ee^{- \omega \VV_\lambda}} = F_\lambda(\omega) + \frac{G_\lambda(\omega)}{1 + H_\lambda(\omega)}
\hspace*{.4in} \mbox{ Re } \omega \geq 0,
\label{form1}
\end{equation}
for any $\lambda < \lambda^*$ with $F_\lambda(0) + G_\lambda(0) = 1$ and $H_\lambda(0) = 0$, for suitably defined functions $F_\lambda(\cdot)$, $G_\lambda(\cdot)$ and $H_\lambda(\cdot)$, which may either be explicitly specified or characterized in an implicit manner.

For example, the above form in particular covers the M/G/1 queue with arrival rate~$\lambda$ and service requirement distribution with mean $\beta < \infty$ and LST $\beta(\omega)$.
Thus the traffic intensity is $\rho = \lambda \beta$ and the critical arrival rate is $\lambda^* = 1 / \beta$.
Denote by $\beta^r(\omega) = \frac{1 - \beta(\omega)}{\beta \omega}$ the LST of the residual service requirement distribution.
The Pollaczek-Khinchine formula then yields for any $\rho < 1$, i.e., $\lambda < \lambda^*$,
\begin{eqnarray*}
\expect{\ee^{- \omega \VV_\lambda}}
&=&
\frac{(1 - \rho) \omega}{\omega - \lambda (1 - \beta(\omega))} \\
&=&
%\frac{1 - \rho}{1 - \rho \frac{1 - \beta(\omega)}{\beta \omega}} \\
%&=&
\frac{1 - \rho}{1 - \rho \beta^r(\omega)} \\
&=&
\frac{1 - \rho}{1 - \rho + \rho [1 - \beta^r(\omega)]} \\
&=&
\frac{1}{1 + \frac{\rho}{1 - \rho} [1 - \beta^r(\omega)]} \\
&=&
\frac{1}{1 + \frac{\lambda}{\lambda^* - \lambda} [1 - \beta^r(\omega)]},
\end{eqnarray*}
which matches the form of~(\ref{form1}) with $F_\lambda(\omega) \equiv 0$, $G_\lambda(\omega) \equiv 1$ and $H_\lambda(\omega) = \frac{\lambda}{\lambda^* - \lambda} [1 - \beta^r(\omega)]$. \\

Returning to the general setting, the next proposition will play an instrumental role in identifying a close connection in terms of the representation in~(\ref{form1}) for $\expect{\ee^{- \omega \VV_\lambda}}$ between the tail asymptotics and heavy-traffic limit of~$\VV_\lambda$ in infinite-variance scenarios.

\begin{proposition}

\label{main1}

(i) If, for some given arrival rate $\lambda < \lambda^*$, (a) $F_\lambda(0) - F_\lambda(\omega) = \theta_\lambda \omega^\alpha + \oo(\omega^\alpha)$ as $\omega \downarrow 0$, (b) $G_\lambda(0) - G_\lambda(\omega) = \gamma_\lambda \omega^\alpha + \oo(\omega^\alpha)$ as $\omega \downarrow 0$, and (c) $H_\lambda(\omega) = \kappa_\lambda \omega^\alpha + \oo(\omega^\alpha)$ as $\omega \downarrow 0$, with $\theta_\lambda, \gamma_\lambda, \kappa_\lambda$ finite constants, $\alpha \in (0, 1)$ and $C_{V_\lambda} = \theta_\lambda + \gamma_\lambda + \kappa_\lambda G_\lambda(0) > 0$, then
\begin{equation}
\pr{\VV_\lambda > x} \sim \frac{C_{V_\lambda}}{\Gamma(1 - \alpha)} x^{- \alpha}
\hspace*{.4in} \mbox{ as } x \to \infty.
\label{tailvlx}
\end{equation}
Also, $\VV_\lambda$ has the same tail asymptotics as the stationary workload $\VV_{\hat\lambda}^c$ defined in Section~\ref{modelsetup} with a regularly varying service requirement distribution as in~(\ref{sere1}), when $\nu = \alpha + 1 \in (1, 2)$ and
\[
C_{V_\lambda} = \frac{\hat\lambda C_{\hat{B}}}{c - \hat\lambda \hat\beta}.
\]

%If there exists some $\gamma \in (0, \infty)$ and $\alpha \in (0, 1)$ such that, for every $\omega \geq 0$, (a) $F_\lambda\left(\left(\frac{\lambda^* - \lambda}{\gamma \lambda^*}\right)^{\frac{1}{\alpha}} \omega\right) \to 0$ as $\lambda \uparrow \lambda^*$, (b) $G_\lambda\left(\left(\frac{\lambda^* - \lambda}{\gamma \lambda^*}\right)^{\frac{1}{\alpha}} \omega\right) \to 1$ as $\lambda \uparrow \lambda^*$ and (c) $H_\lambda\left(\left(\frac{\lambda^* - \lambda}{\gamma \lambda^*}\right)^{\frac{1}{\alpha}} \omega\right) \to \omega^\alpha$ as $\lambda \uparrow \lambda^*$, then

(ii) If there exist functions $R_\lambda$ and $S_\lambda$ such that $R_\lambda / S_\lambda \to 1$ as $\lambda \uparrow \lambda^*$ and, for every $\omega \geq 0$, (a) $F_\lambda(R_\lambda \omega) \to 0$, (b) $G_\lambda(R_\lambda \omega) \to 1$, and (c) $H_\lambda(R_\lambda \omega) \to \omega^\alpha$ with $\alpha \in (0, 1)$ as $\lambda \uparrow \lambda^*$, then
\[
S_\lambda \VV_\lambda \to \WW
\hspace*{.4in} \mbox{ as } \lambda \uparrow \lambda^*,
\]
where $\WW$ is a random variable with a Mittag-Leffler distribution (occasionally also referred to as Kovalenko distribution, see for instance~\cite{GK96}) with parameter~$\alpha$, i.e.,
\[
\expect{\ee^{- \omega \WW}} = \frac{1}{1 + \omega^\alpha}
\hspace*{.4in} \mbox{ Re } \omega \geq 0.
\]
Also, if additionally
\[
%\frac{\lambda^* S_\lambda^\alpha}{\lambda^* - \lambda} \frac{C_{\hat{B}}}{\hat\beta} \mbox{ or }
R_\lambda^\alpha \frac{\lambda}{\lambda^* - \lambda} \frac{C_{\hat{B}}}{\hat\beta} \to \zeta
\hspace*{.4in} \mbox{ as } \lambda \uparrow \lambda^*,
\]
for some constant $\zeta \in (0, \infty)$, then $\VV_\lambda$ has the same heavy-traffic limit for $\lambda \uparrow \lambda^*$, up to a relative capacity slack factor~$\zeta$, as $\VV_{\hat\lambda}^c$ defined in Section~\ref{modelsetup} with a regularly varying service requirement distribution as in~(\ref{sere1}) for $\hat\lambda \uparrow \hat\lambda^* = c / \hat\beta$, when $\nu = \alpha + 1 \in (1, 2)$.
Specifically,
\[
\lim_{\epsilon \downarrow 0} \pr{\epsilon \VV_{\lambda(\epsilon)} \leq x} =
\lim_{\epsilon \downarrow 0} \pr{\epsilon \VV_{\hat\lambda(\epsilon)}^c \leq x} =\pr{\WW \leq x}
\hspace*{.4in} \mbox{ as } \epsilon \downarrow 0,
\]
with $\lambda(\epsilon) = \lambda^* \left(1 - \frac{1}{\zeta} \frac{C_{\hat{B}}}{\hat\beta} \epsilon^{\nu - 1}\right)$ and $\hat\lambda(\epsilon) = \hat\lambda^* \left(1 - \frac{C_{\hat{B}}}{\hat\beta} \epsilon^{\nu - 1} \right)$.

%If there exists some $\gamma \in (0, \infty)$ such that, for every $\omega \geq 0$, (a) $F_\lambda\left(\frac{\lambda^* - \lambda}{\gamma \lambda^*} \omega\right) \to 0$ as $\lambda \uparrow \lambda^*$, (b) $G_\lambda\left(\frac{\lambda^* - \lambda}{\gamma \lambda^*} \omega\right) \to 1$ as $\lambda \uparrow \lambda^*$ and (c) $H_\lambda\left(\frac{\lambda^* - \lambda}{\gamma \lambda^*} \omega\right) \to \omega$ as $\lambda \uparrow \lambda^*$, then

(iii) If there exist functions $R_\lambda$ and $S_\lambda$ such that $R_\lambda / S_\lambda \to 1$ as $\lambda \uparrow \lambda^*$ and, for every $\omega \geq 0$, (a) $F_\lambda(R_\lambda \omega) \to 0$, (b) $G_\lambda(R_\lambda \omega) \to 1$, and (c) $H_\lambda(R_\lambda \omega) \to \omega$ as $\lambda \uparrow \lambda^*$, then
\[
S_\lambda \VV_\lambda \to {\bf Y}
\hspace*{.4in} \mbox{ as } \lambda \uparrow \lambda^*,
\]
where ${\bf Y}$ is a unit-mean exponential random variable.

\end{proposition}

\textbf{Proof}

(i) The stated assumption implies that
\begin{eqnarray*}
\expect{\ee^{- \omega \VV_\lambda}} - 1
&=&
F_\lambda(0) - \theta_\lambda \omega^\alpha + (G_\lambda(0) - \gamma_\lambda \omega^\alpha + \oo(\omega^\alpha)) (1 - \kappa_\lambda \omega^\alpha - \oo(\omega^\alpha)) - 1 + \oo(\omega^\alpha) \\
&=&
- (C_{V_\lambda} + \oo(1)) \omega^\alpha
\hspace*{.4in} \mbox{ as } \omega \downarrow 0.
\end{eqnarray*}
The first statement is then obtained by applying Lemma~\ref{tauber} with $n = 0$.
The second statement follows from comparison of~(\ref{tailvlx}) with~(\ref{tailvhlcx},\ref{tailbrx}).

%\expect{\ee^{- \left(\frac{\lambda^* - \lambda}{\gamma \lambda^*}\right)^{\frac{1}{\alpha}} \omega \VV_\lambda}} =
%F_\lambda\left(\left(\frac{\lambda^* - \lambda}{\gamma \lambda^*}\right)^{\frac{1}{\alpha}} \omega\right) +
%\frac{G_\lambda\left(\left(\frac{\lambda^* - \lambda}{\gamma \lambda^*}\right)^{\frac{1}{\alpha}} \omega\right)}{1 + H_\lambda\left(\left(\frac{\lambda^* - \lambda}{\gamma \lambda^*}\right)^{\frac{1}{\alpha}} \omega\right)} \to \frac{1}{1 + \omega^\alpha}
%\hspace*{.4in} \mbox{ as } \lambda \uparrow \lambda^*.

(ii) We obtain
\[
\expect{\ee^{- \omega R_\lambda \VV_\lambda}} =
F_\lambda(R_\lambda \omega) + \frac{G_\lambda(R_\lambda \omega)}{1 + H_\lambda(R_\lambda \omega)} \to \frac{1}{1 + \omega^\alpha}
\hspace*{.4in} \mbox{ as } \lambda \uparrow \lambda^*.
\]
Recognizing the latter expression as the LST of a Mittag-Leffler distribution with parameter~$\alpha$ and noting that $S_\lambda \VV_\lambda = R_\lambda \VV_\lambda S_\lambda / R_\lambda \to R_\lambda \VV_\lambda$ as $\lambda \uparrow \lambda^*$, the first statement then follows from Feller's convergence theorem.
%The second statement is obtained by recalling that the expression for $\expect{\ee^{- \omega \VV_{\hat\lambda}^c}}$ as provided in~(\ref{lstvhlc}) matches the form of~(\ref{form1}) with $F_{\hat\lambda}(\omega) \equiv 0$, $G_{\hat\lambda}(\omega) \equiv 1$ and $H_{\hat\lambda}(\omega) = \frac{\hat\lambda}{\hat\lambda^* - \hat\lambda} [1 - \hat\beta^r(\omega)]$, and noting that
%\beta(\omega) - 1 + \beta \omega = \kappa \omega^\nu + \oo(\omega^\nu)
%\hspace*{.4in} \mbox{ as } \omega \downarrow 0

In order to prove the second statement, recall from Section~\ref{modelsetup} that when we consider the amount of work $\UU_{\hat\lambda}^c = \VV_{\hat\lambda}^c / c$ measured in time units, this system can be viewed as an M/G/1 queue where the service requirements are scaled by a factor~$1 / c$.
Thus, the Pollaczek-Khinchine formula yields, for any $\hat\lambda < \hat\lambda^* = c / \hat\beta$,
\[
\expect{\ee^{- \omega \UU_{\hat\lambda}^c}} =
\frac{1}{1 + \frac{\hat\lambda \hat\beta / c}{1 - \hat\lambda \hat\beta / c} [1 - \beta^{c, r}(\omega)]} =
\frac{1}{1 + \frac{\hat\lambda}{\hat\lambda^* - \hat\lambda} [1 - \beta^{c,r}(\omega)]},
\]
with $\beta^{c,r}(\omega) = \hat\beta^r(\omega / c)$, and
\begin{equation}
\expect{\ee^{- \omega \VV_{\hat\lambda}^c}} = \expect{\ee^{- \omega c \UU_{\hat\lambda}^c}} =
\frac{1}{1 + \frac{\hat\lambda}{\hat\lambda^* - \hat\lambda} [1 - \beta^{c,r}(c \omega)]} =
\frac{1}{1 + \frac{\hat\lambda}{\hat\lambda^* - \hat\lambda} [1 - \hat\beta^r(\omega)]},
\label{lstvhlc}
\end{equation}
which matches the form of~(\ref{form1}) with $F_{\hat\lambda}(\omega) \equiv 0$, $G_{\hat\lambda}(\omega) \equiv 1$ and $H_{\hat\lambda}(\omega) = \frac{\hat\lambda}{\hat\lambda^* - \hat\lambda} [1 - \hat\beta^r(\omega)]$.

As described in Section~\ref{modelsetup}, we have
\[
\hat\beta^r(\omega) - 1 = - \frac{C_{\hat{B}}}{\hat\beta} \omega^{\nu - 1} + \oo(\omega^{\nu - 1})
\hspace*{.4in} \mbox{ as } \omega \downarrow 0,
\]
so that
\begin{eqnarray*}
& &
H_{\hat\lambda}\left(\left(\frac{\hat\lambda^* - \hat\lambda}{\hat\lambda} \frac{\hat\beta}{C_{\hat{B}}}\right)^{\frac{1}{\nu - 1}} \omega\right) \\
&=&
\frac{\hat\lambda}{\hat\lambda^* - \hat\lambda} \left[1 - \hat\beta^r\left(\left(\frac{\hat\lambda^* - \hat\lambda}{\hat\lambda} \frac{\hat\beta}{C_{\hat{B}}}\right)^{\frac{1}{\nu - 1}} \omega\right)\right] \\
&=&
\frac{\hat\lambda}{\hat\lambda^* - \hat\lambda} \left[\frac{C_{\hat{B}}}{\hat\beta} \left(\left(\frac{\hat\lambda^* - \hat\lambda}{\hat\lambda} \frac{\hat\beta}{C_{\hat{B}}}\right)^{\frac{1}{\nu - 1}} \omega\right)^{\nu - 1} +
\oo\left(\left(\left(\frac{\hat\lambda^* - \hat\lambda}{\hat\lambda} \frac{\hat\beta}{C_{\hat{B}}}\right)^{\frac{1}{\nu - 1}} \omega\right)^{\nu - 1}\right)\right] \\
&=&
\omega^{\nu - 1} + \oo(\hat\lambda^* - \hat\lambda) \\
&\to&
\omega^{\nu - 1}
\hspace*{.4in} \mbox{ as } \hat\lambda \uparrow \hat\lambda^*
\end{eqnarray*}
for every $\omega \geq 0$.
It then follows from the first statement that
\[
\left(\frac{\hat\lambda^* - \hat\lambda}{\hat\lambda^*} \frac{\hat\beta}{C_{\hat{B}}}\right)^{\frac{1}{\nu - 1}} \VV_{\hat\lambda}^c
%\left(\frac{\hat\lambda}{\hat\lambda^*} \frac{\hat\lambda^* - \hat\lambda}{\hat\lambda} \frac{\hat\beta}{C_{\hat{B}}}\right)^{\frac{1}{\nu - 1}} \VV_{\hat\lambda}^c =
%\left(\frac{\hat\lambda}{\hat\lambda^*}\right)^{\frac{1}{\nu - 1}} \left(\frac{\hat\lambda^* - \hat\lambda}{\hat\lambda} \frac{\hat\beta}{C_{\hat{B}}}\right)^{\frac{1}{\nu - 1}} \VV_{\hat\lambda}^c
\to \WW
\hspace*{.4in} \mbox{ as } \hat\lambda \uparrow \hat\lambda^*,
\]
and hence
\[
\epsilon \VV_{\hat\lambda(\epsilon)}^c =
\left(\frac{\hat\lambda^* - \hat\lambda(\epsilon)}{\hat\lambda^*} \frac{\hat\beta}{C_{\hat{B}}}\right)^{\frac{1}{\nu - 1}} \VV_{\hat\lambda(\epsilon)}^c \to \WW
\hspace*{.4in} \mbox{ as } \epsilon \downarrow 0.
\]
Since
\[
R_\lambda \left(\frac{\lambda^* - \lambda}{\lambda^*} \zeta \frac{\hat\beta}{C_{\hat{B}}}\right)^{- \frac{1}{\alpha}} =
\left(R_\lambda^\alpha \frac{\lambda^*}{\lambda^* - \lambda} \frac{1}{\zeta} \frac{C_{\hat{B}}}{\hat\beta}\right)^{\frac{1}{\alpha}} =
\left(\frac{\lambda^*}{\lambda}\right)^{\frac{1}{\alpha}}
\left(R_\lambda^\alpha \frac{\lambda}{\lambda^* - \lambda} \frac{C_{\hat{B}}}{\hat\beta}\right)^{\frac{1}{\alpha}} \zeta^{- \frac{1}{\alpha}} \to 1
\]
as $\lambda \uparrow \lambda^*$, it further follows that when $\nu = \alpha + 1$,
\[
\left(\frac{\lambda^* - \lambda}{\lambda^*} \zeta \frac{\hat\beta}{C_{\hat{B}}}\right)^{\frac{1}{\nu - 1}} \VV_\lambda \to \WW
\hspace*{.4in} \mbox{ as } \lambda \uparrow \lambda^*,
\]
and hence
\[
\epsilon \VV_{\lambda(\epsilon)} =
\left(\frac{\lambda^* - \lambda(\epsilon)}{\lambda^*} \zeta \frac{\hat\beta}{C_{\hat{B}}}\right)^{\frac{1}{\nu - 1}} \VV_{\lambda(\epsilon)}
%\left(\frac{\lambda(\epsilon)}{\lambda^*} \frac{\lambda^* - \lambda(\epsilon)}{\lambda(\epsilon) S_{\lambda(\epsilon)}^{\nu - 1}} \zeta \frac{\hat\beta}{C_{\hat{B}}}\right)^{\frac{1}{\nu - 1}} S_{\lambda(\epsilon)} \VV_{\lambda(\epsilon)}
\to \WW
\hspace*{.4in} \mbox{ as } \epsilon \downarrow 0.
\]

%\expect{\ee^{- \left(\frac{\lambda^* - \lambda}{\gamma \lambda^*}\right) \omega \VV_\lambda}} =
%F_\lambda\left(\frac{\lambda^* - \lambda}{\gamma \lambda^*} \omega\right) +
%\frac{G_\lambda\left(\frac{\lambda^* - \lambda}{\gamma \lambda^*} \omega\right)}{1 + H_\lambda\left(\frac{\lambda^* - \lambda}{\gamma \lambda^*} \omega\right)} \to \frac{1}{1 + \omega}
%\hspace*{.4in} \mbox{ as } \lambda \uparrow \lambda^*.

(iii) We obtain
\[
\expect{\ee^{- \omega R_\lambda \VV_\lambda}} =
F_\lambda(R_\lambda \omega) + \frac{G_\lambda(R_\lambda \omega)}{1 + H_\lambda(R_\lambda \omega)} \to \frac{1}{1 + \omega}
\hspace*{.4in} \mbox{ as } \lambda \uparrow \lambda^*.
\]
Noting that $S_\lambda \VV_\lambda = R_\lambda \VV_\lambda S_\lambda / R_\lambda \to R_\lambda \VV_\lambda$ as $\lambda \uparrow \lambda^*$, the statement then follows from Feller's convergence theorem. \\

\textbf{Remark}
It is worth observing that the tail asymptotics and heavy-traffic limit as stated in the above proposition also hold in case $\VV_\lambda$ does not represent the workload but the waiting time for example.
This is immediately evident in case of the M/G/1 queue with unit service speed where these two quantities are in fact equivalent because of the PASTA property, but also true when $\VV_\lambda$ is not directly linked to the workload.
In the latter scenario the asymptotic equivalence with $\VV_{\hat\lambda}^c$ falls beyond the scope of the reduced-load equivalence that will be discussed later, although it can still be viewed as a manifestation of the same underlying averaging principles and scaling properties. \\

Now suppose that conditions (i)(a)-(b)-(c) of Proposition~\ref{main1} are satisfied, with $\kappa_\lambda \in (0, \infty)$ for all $\lambda < \lambda^*$, for some $\alpha \in (0, 1)$ or $\alpha = 1$, and introduce the function
\[
\hat{H}^\alpha(\lambda, \omega) = \frac{H_\lambda(\omega)}{\kappa_\lambda \omega^\alpha}.
\]
Note that $\hat{H}^\alpha(\lambda, \omega) \to 1$ as $\omega \downarrow 0$ by definition of~$\kappa_\lambda$.
If
\[
\hat{H}^\alpha(\lambda, \kappa_\lambda^{- \frac{1}{\alpha}} \omega) \to 1
\hspace*{.4in} \mbox{ as } \lambda \uparrow \lambda^*
\]
for every $\omega \geq 0$, then
\[
H_\lambda(\kappa_\lambda^{- \frac{1}{\alpha}} \omega) =
\kappa_\lambda \left(\kappa_\lambda^{- \frac{1}{\alpha}} \omega\right)^\alpha \hat{H}^\alpha(\lambda, \kappa_\lambda^{- \frac{1}{\alpha}} \omega) =
\omega^\alpha \hat{H}^\alpha(\lambda, \kappa_\lambda^{- \frac{1}{\alpha}} \omega) \to \omega^\alpha
\hspace*{.4in} \mbox{ as } \lambda \uparrow \lambda^*,
\]
for every $\omega \geq 0$, so that conditions~(ii)(c) and (iii)(c) of Proposition~\ref{main1} are satisfied with $R_\lambda = \kappa_\lambda^{- 1 / \alpha}$ for $\alpha \in (0, 1)$ and $\alpha = 1$, respectively.
This is in particular the case if $\kappa_\lambda \to \infty$ as $\lambda \uparrow \lambda^*$ and $\hat{H}^\alpha(\lambda, \omega) \equiv \hat{H}^\alpha(\omega)$ for all $\lambda \in [0, \lambda^*)$, with $\hat{H}^\alpha(\omega) \to 1$ as $\omega \downarrow 0$.

If in addition $F_\lambda\left(\kappa_\lambda^{- 1 / \alpha} \omega\right) \to 0$ and $G_\lambda\left(\kappa_\lambda^{- 1 / \alpha} \omega\right) \to 1$ as $\lambda \uparrow \lambda^*$ for every $\omega \geq 0$, then conditions~(ii)(a)-(b) and (iii)(a)-(b) of Proposition~\ref{main1} are satisfied as well for $\alpha \in (0, 1)$ and $\alpha = 1$, respectively.
This is ensured for example when $F_\lambda(0) \to 0$ with $F_\lambda(\omega) - F_\lambda(0) \leq f_\lambda \omega$ for some constant $f_\lambda$ with $f_\lambda \kappa_\lambda^{- 1 / \alpha} \to 0$ as $\lambda \uparrow \lambda^*$, and $G_\lambda(0) \to 1$ with $G_\lambda(\omega) - G_\lambda(0) \leq g_\lambda \omega$ for some constant $g_\lambda$ with $g_\lambda \kappa_\lambda^{- 1 / \alpha} \to 0$ as $\lambda \uparrow \lambda^*$, as is trivially the case when $F_\lambda(\omega) \equiv 0$ and $G_\lambda(\omega) \equiv 1$.

Under the above conditions, regularly varying tail asymptotics of index~$- \alpha$, $\alpha \in (0, 1)$, go hand in hand with a Mittag-Leffler distribution with parameter~$\alpha$ of the heavy-traffic limit, where the respective pre-factor $\kappa_\lambda$ and scaling factor $S_\lambda$ are related as $\kappa_\lambda \sim S_\lambda^{- \alpha}$.
Also, for $\alpha = 1$, an exponential heavy-traffic limit arises with scaling factor $S_\lambda \sim \kappa_\lambda^{-1}$.
The above conditions apply in many situations of interest, as illustrated by the examples below. \\

%This is in particular the case if $\hat{H}^\alpha(\lambda, \omega) = \hat{H}_1^\alpha(\lambda) \hat{H}_2^\alpha(\omega)$, the function $\hat{H}_1^\alpha(\lambda)$ does not depend on~$\omega$ and is continuous in~$\lambda^*$ with $\hat{H}_1^\alpha(\lambda) \in (0, \infty)$ for all $\lambda \in [0, \lambda^*]$ , and the function $\hat{H}_2^\alpha(\omega)$ does not depend on~$\lambda$ and is continuous in~$0$, with $\hat{H}_2^\alpha(0) \in (0, \infty$.

\textbf{Remark}
Suppose again that conditions (i)(a)-(b)-(c) of Proposition~\ref{main1} are satisfied, with $\kappa_\lambda \in (0, \infty)$ for all $\lambda < \lambda^*$, $F_\lambda\left(\kappa_\lambda^{- 1 / \alpha} \omega\right) \to 0$, $G_\lambda\left(\kappa_\lambda^{- 1 / \alpha} \omega\right) \to 1$ and $H_\lambda\left(\kappa_\lambda^{- 1 / \alpha} \omega\right) \to \omega^\alpha$ as $\lambda \uparrow \lambda^*$ for some $\alpha \in (0, 1)$ for every $\omega \geq 0$.
As noted above, conditions (ii)(a)-(b)-(c) of Proposition~\ref{main1} are then satisfied as well with $S_\lambda = R_\lambda = \kappa_\lambda^{- 1 / \alpha}$.
This further means that $\theta_\lambda \left(\kappa_\lambda^{- 1 / \alpha} \omega\right)^\alpha = \theta_\lambda \omega^\alpha / \kappa_\lambda \to 0$ and $\gamma_\lambda \left(\kappa_\lambda^{- 1 / \alpha} \omega\right)^\alpha = \gamma_\lambda \omega^\alpha / \kappa_\lambda \to 0$ as $\lambda \uparrow \lambda^*$ for every $\omega \geq 0$, and hence $C_{V_\lambda} / \kappa_\lambda \to 1$ as $\lambda \uparrow \lambda^*$.
Statement~(i) yields for any fixed $\lambda < \lambda^*$
\[
\pr{S_\lambda \VV_\lambda > x} = \pr{\VV_\lambda > x / S_\lambda} \sim
\frac{C_{V_\lambda}}{\Gamma(1 - \alpha)} (x / S_\lambda)^{- \alpha}
\hspace*{.4in} \mbox{ as } x \to \infty,
\]
with $C_{V_\lambda} / S_\lambda^{- \alpha} = C_{V_\lambda} / \kappa_\lambda \to 1$ as $\lambda \uparrow \lambda^*$.
Thus,
\[
\lim_{\lambda \uparrow \lambda^*} \lim_{x \to \infty} \pr{S_\lambda \VV_\lambda > x} \Gamma(1 - \alpha) x^\alpha = 1.
\]
Also, statement~(ii) implies $S_\lambda \VV_\lambda \to \WW$ as $\lambda \uparrow \lambda^*$, with
\[
\expect{\ee^{- \omega \WW}} = \frac{1}{1 + \omega^\alpha}
\hspace*{.4in} \mbox{ Re } \omega \geq 0.
\]
Noting that $\expect{\ee^{- \omega \WW}} - 1 = - \omega^\alpha + \oo(\omega^\alpha)$ and invoking Lemma~\ref{tauber} with $n = 0$, we obtain $\pr{\WW > x} \sim \frac{x^{- \alpha}}{\Gamma(1 - \alpha)}$ as $x \to \infty$.
Thus,
\[
\lim_{x \to \infty} \lim_{\lambda \uparrow \lambda^*} \pr{S_\lambda \VV_\lambda > x} \Gamma(1 - \alpha) x^\alpha = 1.
\]
This demonstrates that the tail asymptotics and heavy-traffic limit of $\VV_\lambda$ commute for $\alpha \in (0, 1)$.
In contrast, for $\alpha = 1$, the two limits do not commute.
Indeed, the case $\alpha = 1$ arises for example when service requirements have finite variance, but possibly regularly varying distributions.
While the tail asymptotics are still regularly varying in such settings, the heavy-traffic limit is then exponential.

\begin{example}

Consider an M/G/1 queue with arrival rate~$\lambda$, generic service requirement $\BB$ with mean $\beta < \infty$ and LST $\beta(\omega)$, and critical arrival rate $\lambda^* = 1 / \beta$.
Denote by $\rho = \lambda \beta$ the traffic intensity and by $\beta^r(\omega) = \frac{1 - \beta(\omega)}{\beta \omega}$ the LST of the residual service requirement distribution.
Recall that the expression for $\expect{\ee^{- \omega \VV_\lambda}}$ as provided by the Pollaczek-Khinchine formula matches the form of~(\ref{form1}) with $F_\lambda(\omega) \equiv 0$, $G_\lambda(\omega) \equiv 1$ and $H_\lambda(\omega) = \frac{\lambda}{\lambda^* - \lambda} [1 - \beta^r(\omega)]$. \\

If the service requirement distribution is regularly varying of index~$- \nu$, $\nu \in (1, 2)$, i.e.,
\begin{equation}
\pr{\BB > x} \sim - \frac{C_B}{\Gamma(1 - \nu)} x^{- \nu}
\hspace*{.4in} \mbox{ as } x \to \infty,
\label{sere2}
\end{equation}
then it follows as described in Section~\ref{modelsetup} that
%\beta(\omega) - 1 + \beta \omega = \kappa \omega^\nu + \oo(\omega^\nu)
%\hspace*{.4in} \mbox{ as } \omega \downarrow 0
\[
\beta^r(\omega) - 1 = - \frac{C_B}{\beta} \omega^{\nu - 1} + \oo(\omega^{\nu - 1})
\hspace*{.4in} \mbox{ as } \omega \downarrow 0,
\]
so that conditions~(i)(a)-(b)-(c) of Proposition~\ref{main1} are satisfied with $\alpha = \nu - 1 \in (0, 1)$, $\theta_\lambda = 0$, $\gamma_\lambda = 0$, and
\[
\kappa_\lambda = \frac{\lambda}{\lambda^* - \lambda} \frac{C_B}{\beta} = \frac{\rho}{1 - \rho} \frac{C_B}{\beta}.
\]
Statement~(i) of Proposition~\ref{main1} then yields, for any $\lambda < \lambda^*$,
\[
\pr{\VV_\lambda > x} \sim \frac{\kappa_\lambda}{\Gamma(2 - \nu)} x^{1 - \nu} =
\frac{\rho}{1 - \rho} \frac{C_B}{\Gamma(2 - \nu) \beta} x^{1 - \nu} \sim
\frac{\rho}{1 - \rho} \pr{\BB^r > x}
\hspace*{.4in} \mbox{ as } x \to \infty,
\]
recovering the tail asymptotics obtained in Cohen~\cite{Cohen73}.

Also, $\hat{H}^{\nu - 1}(\lambda, \omega) \equiv \hat{H}^{\nu - 1}(\omega) = [1 - \beta^r(\omega)] \frac{\beta}{C_B} \omega^{1 - \nu} \to 1$ as $\omega \downarrow 0$, ensuring that $H_\lambda\left(\kappa_\lambda^{- \frac{1}{\nu - 1}} \omega\right) \to \omega^{\nu - 1}$ as $\lambda \uparrow \lambda^*$ for every $\omega \geq 0$.
Noting that $\kappa_\lambda \frac{\lambda^* - \lambda}{\lambda^*} \frac{\beta}{C_B} = \frac{\lambda}{\lambda^*} \to 1$ as $\lambda \uparrow \lambda^*$, statement~(ii) of Proposition~\ref{main1} then yields
\[
\left((1 - \rho) \frac{\beta}{C_B}\right)^{\frac{1}{\nu - 1}} \VV_\lambda =
%\left(\frac{\lambda}{\lambda^*} \frac{\lambda^* - \lambda}{\lambda} \frac{\beta}{C_B}\right)^{\frac{1}{\nu - 1}} \VV_\lambda =
%\left(\frac{\lambda}{\lambda^*}\right)^{\frac{1}{\nu - 1}} \kappa_\lambda^{- \frac{1}{\nu - 1}} \VV_\lambda
\left(\frac{\lambda^* - \lambda}{\lambda^*} \frac{\beta}{C_B}\right)^{\frac{1}{\nu - 1}} \VV_\lambda \to \WW
\hspace*{.4in} \mbox{ as } \lambda \uparrow \lambda^*,
\]
where $\WW$ is a random variable with a Mittag-Leffler distribution with parameter $\nu - 1$, which corresponds to the heavy-traffic limit established in Cohen \cite{Cohen97c,Cohen98a} and Boxma \& Cohen \cite{BC99a,BC99b}.

In contrast, when the service requirement distribution has finite second moment $\beta^{(2)} < \infty$,
\[
\beta(\omega) - 1 + \beta \omega = \frac{1}{2} \beta^{(2)} \omega^2 + \oo(\omega^2)
\hspace*{.4in} \mbox{ as } \omega \downarrow 0,
\]
so that
\[
\beta^r(\omega) - 1 = - \frac{\beta^{(2)}}{2 \beta} \omega + \oo(\omega)
\hspace*{.4in} \mbox{ as } \omega \downarrow 0.
\]

Thus $\hat{H}^1(\lambda, \omega) \equiv \hat{H}^1(\omega) = [1 - \beta^r(\omega)] \frac{2 \beta}{\beta^{(2)}} / \omega \to 1$ as $\omega \downarrow 0$, ensuring that $H_\lambda\left(\frac{\lambda^* - \lambda}{\lambda} \frac{\beta^{(2)}}{2 \beta} \omega\right) \to \omega$ as $\lambda \uparrow \lambda^*$ for every $\omega \geq 0$.
Statement (iii) of Proposition~\ref{main1} then yields
\[
(1 - \rho) \frac{2 \beta}{\beta^{(2)}} \VV_\lambda =
%\rho \frac{1 - \rho}{\rho} \frac{2 \beta}{\beta^{(2)}} \VV_\lambda =
%\frac{\lambda}{\lambda^*} \frac{\lambda^* - \lambda}{\lambda} \frac{2 \beta}{\beta^{(2)}} \VV_\lambda
\frac{\lambda^* - \lambda}{\lambda^*} \frac{2 \beta}{\beta^{(2)}} \VV_\lambda \to {\bf Y}
\hspace*{.4in} \mbox{ as } \lambda \uparrow \lambda^*,
\]
where ${\bf Y}$ is a unit-mean exponential random variable, which corresponds to a special case of the heavy-traffic limit for the G/G/1 queue established in the seminal work of Kingman \cite{Kingman61,Kingman62}.
We refer to Kingman~\cite{Kingman65}, Whitt~\cite{Whitt02} and~\cite{Williams16} for a broader overview of heavy-traffic results for queues with service requirement distributions with finite second moment, and will not further consider these cases in the present paper.

\end{example}

\begin{example}

Consider a single-server queue with service speed $c > 0$, a Poisson arrival process of rate~$\hat\lambda$,  and service requirements that are i.i.d.\ copies of a non-negative random variable $\hat{\BB}$ with mean $\hat\beta < \infty$ and LST $\hat\beta(\omega)$, and critical arrival rate $\hat\lambda^* = c / \hat\beta$.
Recall that the expression for $\expect{\ee^{- \omega \VV_{\hat\lambda}^c}}$ as provided in~(\ref{lstvhlc}) matches the form of~(\ref{form1}) with $F_{\hat\lambda}(\omega) \equiv 0$, $G_{\hat\lambda}(\omega) \equiv 1$ and $H_{\hat\lambda}(\omega) = \frac{\hat\lambda}{\hat\lambda^* - \hat\lambda} [1 - \hat\beta^r(\omega)]$. \\

If the service requirement distribution is regularly varying of tail index~$- \nu$, $\nu \in (1, 2)$, as in~(\ref{sere1}), then it follows as described in Section~\ref{modelsetup} that
%\beta(\omega) - 1 + \beta \omega = \kappa \omega^\nu + \oo(\omega^\nu)
%\hspace*{.4in} \mbox{ as } \omega \downarrow 0,
\[
\hat\beta^r(\omega) - 1 = - \frac{C_{\hat{B}}}{\hat\beta} \omega^{\nu - 1} + \oo(\omega^{\nu - 1})
\hspace*{.4in} \mbox{ as } \omega \downarrow 0,
\]
so that conditions~(i)(a)-(b)-(c) of Proposition~\ref{main1} are satisfied with $\alpha = \nu - 1 \in (0, 1)$, $\theta_{\hat\lambda} = 0$, $\gamma_{\hat\lambda} = 0$, and
\[
\kappa_{\hat\lambda} = \frac{\hat\lambda}{\hat\lambda^* - \hat\lambda} \frac{C_{\hat{B}}}{\hat\beta} = \frac{\hat\lambda \hat\beta}{c - \hat\lambda \hat\beta} \frac{C_{\hat{B}}}{\hat\beta}.
\]
Statement~(i) of Proposition~\ref{main1} then yields, for any $\hat\lambda < \hat\lambda^*$,
\[
\pr{\VV_{\hat\lambda} > x} \sim \frac{\kappa_{\hat\lambda}}{\Gamma(2 - \nu)} x^{1 - \nu} =
\frac{\hat\lambda \hat\beta}{c - \hat\lambda \hat\beta} \frac{C_{\hat{B}}}{\Gamma(2 - \nu) \hat\beta} x^{1 - \nu}
\hspace*{.4in} \mbox{ as } x \to \infty,
\]
recovering~(\ref{tailvhlcx}).

Also, $\hat{H}^{\nu - 1}(\hat\lambda, \omega) \equiv \hat{H}^{\nu - 1}(\omega) = [1 - \beta^r(\omega)] \frac{C_{\hat{B}}}{\hat\beta} \omega^{1 - \nu} \to 1$ as $\omega \downarrow 0$, ensuring that $H_{\hat\lambda}\left(\kappa_{\hat\lambda}^{- \frac{1}{\nu - 1}} \omega\right) \to \omega^{\nu - 1}$ as $\hat\lambda \uparrow \hat\lambda^*$ for every $\omega \geq 0$.
Noting that $\kappa_{\hat\lambda} \frac{\hat\lambda^* - \hat\lambda}{\hat\lambda^*} \frac{\hat\beta}{C_{\hat{B}}} = \frac{\hat\lambda}{\hat\lambda^*} \to 1$ as $\hat\lambda \uparrow \hat\lambda^*$, statement~(ii) of Proposition~\ref{main1} then yields
\[
\left(\left(1 - \frac{\hat\lambda \hat\beta}{c}\right) \frac{\hat\beta}{C_{\hat{B}}}\right)^{\frac{1}{\nu - 1}} \VV_\lambda =
%\left(\frac{\hat\lambda \beta}{c} \frac{c - \hat\lambda \beta}{\hat\lambda \beta} \frac{\beta}{C_B}\right)^{\frac{1}{\nu - 1}} \VV_\lambda =
%\left(\frac{\hat\lambda}{\hat\lambda^*} \frac{\hat\lambda^* - \hat\lambda}{\hat\lambda} \frac{\hat\beta}{C_{\hat{B}}}\right)^{\frac{1}{\nu - 1}} \VV_\lambda =
%\left(\frac{\hat\lambda}{\hat\lambda^*}\right)^{\frac{1}{\nu - 1}} \hat\kappa_{\hat\lambda}^{- \frac{1}{\nu - 1}} \VV_\lambda
\left(\frac{\hat\lambda^* - \hat\lambda}{\hat\lambda^*} \frac{\hat\beta}{C_{\hat{B}}}\right)^{\frac{1}{\nu - 1}} \VV_\lambda \to \WW
\hspace*{.4in} \mbox{ as } \lambda \uparrow \lambda^*,
\]
where $\WW$ is a random variable with a Mittag-Leffler distribution with parameter $\nu - 1$.

\end{example}

\begin{example}

Next consider a multi-class M/G/1 queue with $K$~customer classes and total arrival rate~$\lambda$.
The fraction of class-$i$ customers is $p_i$, and the generic class-$i$ service requirement is a non-negative random variable $\BB_i$ with mean $\beta_i < \infty$ and LST $\beta_i(\omega)$, $i = 1, \dots, K$.
Denote by $\rho_i = \lambda p_i \beta_i$ the traffic intensity of class~$i$, and by $\rho = \sum_{i = 1}^{K} \rho_i = \lambda \beta$ the total traffic intensity, with $\beta = \sum_{i = 1}^{K} p_i \beta_i$ the mean service requirement of an arbitrary customer, so that the critical arrival rate is $\lambda^* = 1 / \beta$.
Denote the LST of the service requirement distribution of an arbitrary customer by $\beta(\omega) = \sum_{i = 1}^{K} p_i \beta_i(\omega)$, and that of the residual service requirement distribution by
\[
\beta^r(\omega) = \frac{1 - \beta(\omega)}{\beta \omega} =
\frac{1 - \sum_{i = 1}^{K} p_i \beta_i(\omega)}{\sum_{i = 1}^{K} p_i \beta_i \omega} =
\frac{\sum_{i = 1}^{K} p_i (1 - \beta_i(\omega))}{\sum_{i = 1}^{K} p_i \beta_i \omega}.
\]
As before, the expression for $\expect{\ee^{- \omega \VV_\lambda}}$ matches the form of~(\ref{form1}) with $F_\lambda(\omega) \equiv 0$, $G_\lambda(\omega) \equiv 1$ and $H_\lambda(\omega) = \frac{\lambda}{\lambda^* - \lambda} [1 - \beta^r(\omega)]$. \\

Now suppose that for one of the classes~$i_0$ the service requirement distribution is regulary varying of index~$- \nu$, $\nu \in (1, 2)$, i.e.,
\[
\pr{\BB_{i_0} > x} \sim - \frac{C_0}{\Gamma(1 - \nu)} x^{- \nu}
\hspace*{.4in} x \to \infty,
\]
while $\pr{\BB_i > x} = \oo(x^{- \nu})$ as $x \to \infty$ for all classes $i \neq i_0$.
It then follows from Lemma~\ref{tauber} with $\alpha = \nu$ and $n = 1$ that
\[
\beta(\omega) - 1 + \beta \omega = \alpha_{i_0} C_0 \omega^\nu + \oo(\omega^\nu)
\hspace*{.4in} \mbox{ as } \omega \downarrow 0,
\]
and thus
\[
\beta^r(\omega) - 1 = - \frac{\alpha_{i_0} C_0}{\beta} \omega^{\nu - 1} + \oo(\omega^{\nu - 1})
\hspace*{.4in} \mbox{ as } \omega \downarrow 0,
\]
so that conditions~(i)(a)-(b)-(c) of Proposition~\ref{main1} are satisfied with $\alpha = \nu - 1 \in (0, 1)$, $\theta_\lambda = 0$, $\gamma_\lambda = 0$, and
\[
\kappa_\lambda = \frac{\lambda}{\lambda^* - \lambda} \frac{p_{i_0} C_0}{\beta} =
\frac{\rho}{1 - \rho} \frac{p_{i_0} C_0}{\beta} =
\frac{\lambda p_{i_0} C_0}{1 - \rho} = \frac{\rho_{i_0}}{1 - \rho} \frac{C_0}{\beta_{i_0}}.
\]
Statement~(i) of Proposition~\ref{main1} then yields that
\[
\pr{\VV_\lambda > x} \sim \frac{\kappa_\lambda}{\Gamma(2 - \nu)} x^{1 - \nu} =
%\frac{\rho}{1 - \rho} \frac{p_{i_0} C_0}{\Gamma(2 - \nu) \beta} x^{1 - \nu} =
\frac{\rho_{i_0}}{1 - \rho} \frac{C_0}{\Gamma(2 - \nu) \beta_{i_0}} x^{1 - \nu}
%\frac{\rho_{i_0}}{c - \rho_{i_0}} \pr{\BB_{i_0}^r > x},
\hspace*{.4in} \mbox{ as } x \to \infty,
\]
and that $\VV_\lambda$ has the same tail asymptotics as the stationary workload $\VV_{\hat\lambda}^c$ defined in Section~\ref{modelsetup} with a regularly varying service requirement distribution of index $- \nu$, $\nu = \alpha + 1 \in (1, 2)$, as in~(\ref{sere1}) when
\[
\kappa_\lambda = %\frac{\rho}{1 - \rho} \frac{p_{i_0} C_0}{\beta} =
\frac{p_{i_0} \lambda C_0}{1 - \rho} = \frac{\hat\lambda C_{\hat{B}}}{c - \hat\lambda \hat\beta}.
\]
The above equation is satisfied when we take $\hat\lambda = \alpha_{i_0} \lambda$, $\hat{\BB}$ to be the service requirement $\BB_{i_0}$ of class-$i_0$ customers so that $\hat\beta = \beta_{i_0}$ and $C_{\hat{B}} = C_0$, and
\[
c = 1 - \rho + \hat\lambda \hat\beta = 1 - \rho + p_{i_0} \lambda \beta_{i_0} =
1 - \rho + \rho_{i_0} = 1 - \sum_{i \neq i_0} \rho_i
\]
representing the full service rate reduced by the aggregate load of all classes $i \neq i_0$.
In other words, the tail asymptotics of the workload are similar to those in a reduced system handling traffic of class~$i_0$ only and having service speed~$c$, which is reminiscent of the so-called reduced-load equivalence established in~\cite{AMN99}.

Also, $\hat{H}^{\nu - 1}(\lambda, \omega) \equiv \hat{H}^{\nu - 1}(\omega) = [1 - \beta^r(\omega)] \frac{\beta}{p_{i_0} C_0} \omega^{1 - \nu} \to 1$ as $\omega \downarrow 0$, ensuring that $H_\lambda\left(\kappa_\lambda^{- \frac{1}{\nu - 1}} \omega\right) \to \omega^{\nu - 1}$ as $\lambda \uparrow \lambda^*$ for every $\omega \geq 0$.
Noting that $\kappa_\lambda \frac{\lambda^* - \lambda}{\lambda^*} \frac{\beta}{p_{i_0} C_0} = \frac{\lambda}{\lambda^*} \to 1$ as $\lambda \uparrow \lambda^*$, statement~(ii) of Proposition~\ref{main1} then yields
\[
\left((1 - \rho) \frac{\beta}{p_{i_0} C_0}\right)^{\frac{1}{\nu - 1}} \VV_\lambda =
%\left(\rho \frac{1 - \rho}{\rho} \frac{\beta}{p_{i_0} C_0}\right)^{\frac{1}{\nu - 1}} \VV_\lambda =
%\left(\frac{\lambda}{\lambda^*} \frac{\lambda^* - \lambda}{\lambda} \frac{\beta}{p_{i_0} C_0}\right)^{\frac{1}{\nu - 1}} \VV_\lambda =
%\left(\frac{\lambda}{\lambda^*}\right)^{\frac{1}{\nu - 1}} \kappa_\lambda^{- \frac{1}{\nu - 1}} \VV_\lambda
\left(\frac{\lambda^* - \lambda}{\lambda^*} \frac{\beta}{p_{i_0} C_0}\right)^{\frac{1}{\nu - 1}} \VV_\lambda \to \WW
\hspace*{.4in} \mbox{ as } \lambda \uparrow \lambda^*,
\]
where $\WW$ is a random variable with a Mittag-Leffler distribution with parameter $\nu - 1$.

Additionally, as
\[
\frac{\lambda}{(\lambda^* - \lambda) \kappa_\lambda} \frac{C_{\hat{B}}}{\hat\beta} = \zeta =
\frac{\beta}{p_{i_0} C_0} \frac{C_{\hat{B}}}{\hat\beta},
\]
we find that $\VV_\lambda$ has the same heavy-traffic limit for $\lambda \uparrow \lambda^*$, up to a relative capacity slack factor~$\zeta$, as $\VV_{\hat\lambda}^c$ defined in Section~\ref{modelsetup} with a regularly varying service requirement distribution of index $- \nu$, $\nu = \alpha + 1 \in (1, 2)$, as in~(\ref{sere1}) for $\hat\lambda \uparrow \hat\lambda^* = c / \hat\beta$.
This holds in particular when
\[
c = c^* = 1 - \sum_{i \neq i_0} \lambda^* p_i \beta_i = \lambda^* \sum_{i = 1}^{K} p_i \beta_i - \sum_{i \neq i_0} \lambda^* p_i \beta_i = \lambda^* p_{i_0} \beta_{i_0},
\]
representing the full service rate reduced by the aggregate load of all classes $i \neq i_0$ in heavy-traffic conditions as $\lambda \uparrow \lambda^*$.
If we take $\hat\lambda = p_{i_0} \lambda$ and $\hat{\BB}$ to be the service requirement $\BB_{i_0}$ of class-$i_0$ customers so that $\hat\beta = \beta_{i_0}$ and $C_{\hat{B}} = C_0$, then we conclude that $\VV_\lambda$ has the same heavy-traffic limit, up to a relative capacity slack factor $\zeta = \frac{\beta}{p_{i_0} \beta_{i_0}}$, as the stationary workload in the reduced system described above with service speed~$c^*$.
This may be interpreted as a heavy-traffic counterpart of the reduced-load equivalence for the tail asymptotics.
Specifically, if we define
\[
\lambda(\epsilon) =
\lambda^* \left(1 - \frac{1}{\zeta} \frac{C_0}{\beta_{i_0}} \epsilon^{\nu - 1}\right) =
%\lambda^* \left(1 - \frac{p_{i_0} \beta_{i_0}}{\beta} \frac{C_0}{\beta_{i_0}} \epsilon^{\nu - 1}\right) =
\lambda^* \left(1 - \frac{p_{i_0} C_0}{\beta} \epsilon^{\nu - 1}\right) =
\lambda^* \frac{\beta - p_{i_0} C_0 \epsilon^{\nu - 1}}{\beta} =
\frac{1 - \lambda^* p_{i_0} C_0 \epsilon^{\nu - 1}}{\beta}
\]
and
\[
\hat\lambda(\epsilon) =
\hat\lambda^* \left(1 - \frac{C_0}{\beta_{i_0}} \epsilon^{\nu - 1}\right) =
p_{i_0} \lambda^* \frac{\beta_{i_0} - C_0 \epsilon^{\nu - 1}}{\beta_{i_0}} =
\frac{c^* - \lambda^* p_{i_0} C_0 \epsilon^{\nu - 1}}{\beta_{i_0}},
\]
then
\[
\lim_{\epsilon \downarrow 0} \pr{\epsilon \VV_{\lambda(\epsilon)} \leq x} =
\lim_{\epsilon \downarrow 0} \pr{\epsilon \VV_{\hat\lambda(\epsilon)}^{c^*} \leq x} = \pr{\WW \leq x}.
\]
Note that $\frac{1}{\zeta} = \frac{p_{i_0} \beta_{i_0}}{\beta}$ represents the heaviest-tail load as fraction of the total load, and that $\lambda(\epsilon)$ and $\hat\lambda(\epsilon)$ correspond to the same absolute capacity slack $\lambda^* p_{i_0} C_0 \epsilon^{\nu - 1}$ with respect to the service speeds~$1$ and~$c^*$ in the two respective systems.

\end{example}

\section{M/G/1 queue with alternating service speed}
\label{alternating}

In this section we consider an M/G/1 queue with a time-varying service speed as studied by Boxma \& Kurkova~\cite{BK01}.
Customers arrive as a Poisson process of rate~$\lambda$ and have a general service requirement distribution $B(x)$ with LST $\beta(\omega) = \expect{\ee^{- \omega \BB}}$, mean~$\beta < \infty$ and second moment $\beta^{(2)}$ which may or may not be infinite.
The service speed toggles between a low value~$s_L$ and a high value~$s_H$ according to an alternating renewal process.
The high-speed periods are exponentially distributed with mean $1/\nu$, while the low-speed periods have a general distribution $D(x)$ with LST $\delta(\omega) = \expect{\ee^{- \omega \DD}}$, mean $\delta < \infty$, and second moment $\delta^{(2)}$ which may or may not be infinite.

Henceforth, the stability condition $\lambda \beta < \bar{s}$ is assumed to be satisfied, where
\[
\bar{s} := \frac{\delta}{\delta + 1/\nu} s_L + \frac{1/\nu}{\delta + 1/\nu} s_H =
\frac{\nu \delta}{1 + \nu \delta} s_L + \frac{1}{1 + \nu \delta} s_H =
\frac{s_H + \nu \delta s_L}{1 + \nu \delta}
\]
represents the time-average service speed, so the critical arrival rate is $\lambda^* = \bar{s} / \beta$.
In addition, we will assume that $\lambda \beta > s_L$, implying that the workload has positive drift during low-speed periods and that the high-speed periods are thus essential for stability. \\
%We will examine the behavior of the workload in a heavy-traffic regime where $\epsilon := \bar{s} - \lambda \beta \downarrow 0$ as the arrival rate~$\lambda$ approaches the critical value $\lambda_0 = \bar{s} / \beta$.

Let $\VV_\lambda$ be the workload in the system and let the binary random variable $\II$ indicate whether the system is in a high-speed period ($\II = H$) or a low-speed period ($\II = L$).
Let $F_H(x) := \pr{\VV_\lambda < x; \II = H}$ be the equilibrium probability that the workload~$\VV_\lambda$ is less than~$x$ and that the service speed is high, and let $\Phi_H(\omega) := \expect{\ee^{- \omega \VV_\lambda}; \II = H} = \int_{x = 0}^{\infty} \ee^{- \omega x} \dd F_H(x)$.
Note that $\Phi_H(0) = \pr{\II = H} = \frac{1/\nu}{\delta + 1/\nu} = \frac{1}{1 + \nu \delta}$.
Further let $\VV_{H, \lambda}$ be a random variable with the conditional distribution of the workload~$\VV_\lambda$ given that the service speed is high, i.e., $\pr{\VV_{H, \lambda} < x} = \pr{\VV_\lambda < x | \II = H} = \pr{\VV_\lambda < x; \II = H} / \pr{\II = H} = (1 + \nu \delta) \pr{\VV_\lambda < x; \II = H}$, so that $\expect{\ee^{- \omega \VV_{H,\lambda}}} = \expect{\ee^{- \omega \VV_\lambda} | \II = H} = \Phi_H(\omega) / \Phi_H(0) = (1 + \nu \delta) \Phi_H(\omega)$. \\

For later reference, we introduce a closely related system with constant service speed $s_H$ which besides the instantaneous input associated with the customers as described above is also fed by gradual input as produced by an On-Off source, which generates traffic at rate $r = s_H - s_L$ when On.
The On- and Off-periods correspond to the low-speed and high-speed periods in the original system, respectively.
Thus the fraction of time that the source is On is $p_{On} = \frac{\delta}{\delta + 1 / \nu} = \frac{\nu \delta}{1 + \nu \delta}$, and the time-average traffic rate is $p_{On} r = \frac{\nu \delta (s_H - s_L)}{1 + \nu \delta}$.

We assume that the traffic from the On-Off source receives priority over the traffic from the customers.
This priority rule does not affect the total workload, but implies that there is never any workload from the On-Off source.
The workload thus entirely consists of the workload from the customers, and is served at rate $s_L$ or $s_H$ when the On-Off source is On or Off, respectively, so that it evolves exactly like the workload in the original system with the time-varying service speed.
In particular, $\VV_\lambda$ may equivalently be interpreted as the workload in the fixed-speed system with additional gradual input, and $\VV_{H, \lambda}$ also represents the workload in the latter system in equilibrium given that the On-Off source is Off.

For conciseness, we will henceforth refer to the original system as the \textit{variable-speed system} and to the corresponding fixed-speed system with additional gradual input as the \textit{dual-input system}.
In addition we introduce two reference systems:
Reference system~A is a system with constant service speed $c = \bar{s}$ fed by the instantaneous input from the customers only;
Reference system~B is a system with constant service speed $d = s_H - \lambda \beta$ fed by the gradual input from the above-described On-Off source only. \\

Equation (2.18) in~\cite{BK01} yields:
\begin{equation}
\Phi_H(\omega) K_\lambda(\omega) = - s_H \omega F_H(0) - s_L \omega \Psi_\lambda(\omega),
\label{funcequa1}
\end{equation}
$\mbox{Re } \omega \geq 0$, where
\[
K_\lambda(\omega) =
\nu + \lambda (1 - \beta(\omega)) - s_H \omega - \nu \delta(\lambda (1 - \beta(\omega)) - s_L \omega),
\]
and the function $\Psi_\lambda(\omega)$ is defined in Equation (2.19) of~\cite{BK01} (be it that somewhat different notation is used there).
Note that
\begin{eqnarray*}
K_\lambda(\omega)
&=& \nu + \lambda [\beta \omega + \oo(\omega)] - s_H \omega - \nu [1 - \delta [\lambda \beta \omega - s_L \omega + \oo(\omega)] + \oo(\omega)]
\hspace*{.4in} \mbox{ as } \omega \downarrow 0 \\
%\lambda \beta \omega - s_H \omega + \nu \delta [\lambda \beta \omega - s_L \omega] + \oo(\omega) \\
&=&
(1 + \nu \delta) \lambda \beta \omega - s_H \omega - \nu \delta s_L \omega + \oo(\omega) \\
&=&
(1 + \nu \delta) (\lambda \beta - \bar{s}) \omega + \oo(\omega).
\end{eqnarray*}

Thus, substituting $\omega = 0$ in~(\ref{funcequa1}) and noting that $\Phi_H(0) = 1 / (1 + \nu \delta)$, we obtain
\[
- s_H F_H(0) - s_L \Psi_\lambda(0) = \lambda \beta - \bar{s},
\]
yielding
\begin{eqnarray*}
\expect{\ee^{- \omega \VV_{H, \lambda}}}
&=&
(1 + \nu \delta) \Phi_H(\omega) \\
&=&
- \frac{(1 + \nu \delta) [s_H \omega F_H(0) + s_L \omega \Psi_\lambda(\omega)]}{K_\lambda(\omega)} \\
&=&
- \frac{(1 + \nu \delta) [s_H F_H(0) + s_L \Psi_\lambda(0) + s_L [\Psi_\lambda(\omega) - \Psi_\lambda(0)]] \omega}{K_\lambda(\omega)} \\
&=&
- \frac{(1 + \nu \delta) [\bar{s} - \lambda \beta + s_L [\Psi_\lambda(\omega) - \Psi_\lambda(0)]] \omega}{K_\lambda(\omega)},
\end{eqnarray*}
which matches the form of~(\ref{form1}), with $F_\lambda(\omega) \equiv 0$,
\[
G_\lambda(\omega) = 1 - \frac{s_L [\Psi_\lambda(0) - \Psi_\lambda(\omega)]}{\bar{s} - \lambda \beta},
\]
and
\[
H_\lambda(\omega) = - \frac{1}{(1 + \nu \delta) (\bar{s} - \lambda \beta)} \left[\frac{K_\lambda(\omega)}{\omega} + (1 + \nu \delta) (\bar{s} - \lambda \beta)\right].
\]

\iffalse

\begin{eqnarray*}
\hat{H}^\alpha(\lambda, \omega)
&=&
\frac{\lambda^* - \lambda}{\lambda \omega^\alpha} H_\lambda(\omega) \\
&=&
\frac{\lambda^* - \lambda}{\lambda \omega^\alpha} \frac{1}{(1 + \nu \delta) (\lambda \beta - \bar{s})} \left[\frac{K_\lambda(\omega)}{\omega} - (1 + \nu \delta) (\lambda \beta - \bar{s})\right] \\
&=&
\frac{\lambda^* - \lambda}{\lambda \omega^\alpha} \frac{1}{(1 + \nu \delta) (\lambda - \bar{s} / \beta) \beta} \left[\frac{K_\lambda(\omega)}{\omega} - (1 + \nu \delta) (\lambda \beta - \bar{s})\right] \\
&=&
\frac{\lambda^* - \lambda}{\lambda \omega^\alpha} \frac{1}{(1 + \nu \delta) (\lambda - \lambda^*) \beta} \left[\frac{K_\lambda(\omega)}{\omega} - (1 + \nu \delta) (\lambda \beta - \bar{s})\right] \\
&=&
\frac{1}{(1 + \nu \delta) \lambda \beta \omega^\alpha} \left[\frac{K_\lambda(\omega)}{\omega} - (1 + \nu \delta) (\lambda \beta - \bar{s})\right].
\end{eqnarray*}

\fi

In~\cite{BK01} it is shown that $\Psi_\lambda(0) - \Psi_\lambda(\omega) = \OO(\omega)$ as $\omega \downarrow 0$ for any $\lambda < \lambda^*$, so that $1 - G_\lambda(\omega) = \oo(\omega^\alpha)$ as $\omega \downarrow 0$ for any $\alpha \in (0, 1)$, and that
\[
\frac{K_\lambda(\omega)}{\omega} + (1 + \nu \delta) (\bar{s} - \lambda \beta) =
- \eta_\lambda \omega^\alpha + \oo(\omega^\alpha)
\hspace*{.4in} \mbox{ as } \omega \downarrow 0,
\]
with \\
(a) $\alpha = \nu_B - 1\in (0, 1)$ and $\eta_\lambda = \lambda (1 + \nu \delta) C_B$ if $1 - B(x) \sim - \frac{C_B}{\Gamma(1 - \nu_B)} x^{- \nu_B}$ with $\nu_B \in (1, 2)$ and $1 - D(x) = \oo(x^{- \nu_B})$ as $x \to \infty$; \\
(b) $\alpha = \nu_D - 1 \in (0, 1)$ and $\eta_\lambda = \nu (\lambda \beta - s_L)^{\nu_D} C_D$ if $1 - D(x) \sim - \frac{C_D}{\Gamma(1 - \nu_D)} x^{- \nu_D}$ with $\nu_D \in (1, 2)$ and $1 - B(x) = \oo(x^{- \nu_D})$ as $x \to \infty$; \\
(c) $\alpha = \nu_0 - 1 \in (0, 1)$ and $\eta_\lambda = \lambda (1 + \nu \delta) C_B + \nu (\lambda \beta - s_L)^{\nu_0} C_D$ if $1 - B(x) \sim - \frac{C_B}{\Gamma(1 - \nu_0)} x^{- \nu_0}$ and $1 - D(x) \sim - \frac{C_D}{\Gamma(1 - \nu_0)} x^{- \nu_0}$ with $\nu_0 \in (1, 2)$ as $x \to \infty$.

Thus $H_\lambda(\omega) = \kappa_\lambda \omega^\alpha + \oo(\omega^\alpha)$ as $\omega \downarrow 0$, with
\[
\kappa_\lambda = \frac{\eta_\lambda}{(1 + \nu \delta) (\bar{s} - \lambda \beta)} =
\frac{\eta_\lambda}{(1 + \nu \delta) (\lambda^* - \lambda) \beta}.
\]

\iffalse

\[
\hat{H}^\alpha(\lambda, \omega) \to \hat\kappa_\lambda
\hspace*{.4in} \mbox{ as } \omega \downarrow 0
\]
for all $\lambda \in [0, \lambda^*]$, and
\[
\hat{H}\left(\lambda, \left(\frac{\lambda^* - \lambda}{\gamma \lambda^*}\right)^{\frac{1}{\alpha}}\right) \to \hat\kappa_{\lambda^*}
\hspace*{.4in} \mbox{ as } \lambda \uparrow \lambda^*
\]
for every $\omega \geq 0$, with
\[
\hat\kappa_\lambda = \frac{D_\lambda}{(1 + \nu \delta) \lambda \beta},
\]
\[
\gamma = \hat\kappa_{\lambda^*} = \frac{D_{\lambda^*}}{(1 + \nu \delta) \lambda^* \beta} = \frac{D_{\lambda^*}}{(1 + \nu \delta) \bar{s}}.
\]

\fi

\subsection{Tail asymptotics}

We first consider the tail asymptotics of the workload as obtained in~\cite{BK01}, but we consider these through the lens of~(\ref{form1}) in order to highlight the connection with the heavy-traffic limit.
Statement~(i) of Proposition~\ref{main1} implies
\[
\pr{\VV_{H, \lambda} > x} \sim \frac{C_{V_{H,\lambda}}}{\Gamma(1 - \alpha)} x^{- \alpha}
\hspace*{.4in} \mbox{ as } x \to \infty,
\]
with $C_{V_{H,\lambda}} = \frac{\eta_\lambda}{(1 + \nu \delta) (\bar{s} - \lambda \beta)}$, which is consistent with Theorem~4.1 in~\cite{BK01}.

Also, $\VV_{H, \lambda}$ has the same tail asymptotics as the stationary workload $\VV_{\hat\lambda}^c$ defined in Section~\ref{modelsetup} with a regularly varying service requirement distribution of index $- \nu$, $\nu = 1 + \alpha \in (1, 2)$, as in~(\ref{sere1}) when
\begin{equation}
C_{V_{H,\lambda}} = \frac{\hat\lambda C_{\hat{B}}}{c - \hat\lambda \hat\beta}.
\label{cvhl1}
\end{equation}

In case~(a), we have $\alpha = \nu_B - 1 \in (0, 1)$ and $\eta_\lambda = \lambda (1 + \nu \delta) C_B$, so that $C_{V_{H,\lambda}} = \frac{\lambda C_B}{\bar{s} - \lambda \beta}$, and
\begin{eqnarray*}
\pr{\VV_{H, \lambda} > x}
&\sim&
\frac{\lambda C_B}{(\bar{s} - \lambda \beta) \Gamma(2 - \nu_B)} x^{1 - \nu_B}
\hspace*{.4in} \mbox{ as } x \to \infty \\
&=&
\frac{\lambda \beta}{\bar{s} - \lambda \beta} \frac{C_B}{\Gamma(2 - \nu_B) \beta} x^{1 - \nu_B}.
\end{eqnarray*}

The equality in~(\ref{cvhl1}) is satisfied when $\hat\lambda = \lambda$, $\hat{\BB} = \BB$, and $c = \bar{s}$, implying that $\VV_{H, \lambda}$ has the same tail asymptotics as $\VV_{\hat\lambda}^c$ in that case, as is confirmed by comparison with~(\ref{tailvhlcx}).
Thus, in case~(a) $\VV_{H, \lambda}$ has the same tail asymptotics as the stationary workload in Reference System~A described above.
Now recall that $\VV_{H, \lambda}$ may also be interpreted as the workload in the dual-input system mentioned earlier, given that the On-Off source in that system is Off, and hence the tail asymptotics of the workload in the latter system are also the same as in Reference System~A.
This property is a manifestation of the reduced-load equivalence for the workload asymptotics in the dual-input system when the instantaneous input is `heavier-tailed' than the gradual input.
Observe here that the service speed~$\bar{s}$ of Reference System~A equals the service speed $s_H$ of the dual-input system reduced by the average rate $\frac{\nu \delta (s_H - s_L)}{1 + \nu \delta}$ of the On-Off source:
\[
s_H - \frac{\nu \delta (s_H - s_L)}{1 + \nu \delta} = \frac{s_H + \nu \delta s_L}{1 + \nu \delta} = \bar{s}.
\]

In case~(b), we have $\alpha = \nu_D - 1 \in (0, 1)$ and $\eta_\lambda = \nu (\lambda \beta - s_L)^{\nu_D} C_D$, so that
\begin{eqnarray*}
\pr{\VV_{H, \lambda} > x}
&\sim&
\frac{\nu (\lambda \beta - s_L)^{\nu_D} C_D}{(1 + \nu \delta) (\bar{s} - \lambda \beta) \Gamma(2 - \nu_D)} x^{1 - \nu_D}
\hspace*{.4in} \mbox{ as } x \to \infty \\
&=&
\frac{\nu \delta}{1 + \nu \delta} \frac{(\lambda \beta - s_L)^{\nu_D}}{\bar{s} - \lambda \beta} \frac{C_D}{\Gamma(2 - \nu_D) \delta} x^{1 - \nu_D}.
\end{eqnarray*}

Now let us consider Reference System~B described above, where the workload increases at rate $r - d = s_H - s_L - (s_H - \lambda \beta) = \lambda \beta - s_L > 0$ during On-periods and decreases at rate $d = s_H - \lambda \beta > 0$ during Off-periods.
If we denote by $\WW_{On}^d$ the workload in Reference System~B at the start of an On-period and assume that $D(x)$ is regularly varying with tail index $- \nu_D$ as in case~(b), then it follows from the results in~\cite{BD98a} that
\[
\pr{\WW_{On}^d > x} \sim
\frac{\nu \delta}{1 + \nu \delta} \frac{(\lambda \beta - s_L)^{\nu_D}}{\bar{s} - \lambda \beta}
\frac{C_D}{\Gamma(2 - \nu_D) \delta} x^{1 - \nu_D}
\hspace*{.4in} \mbox{ as } x \to \infty.
\]
%\begin{eqnarray*}
%\pr{\WW_{On}^d > x}
%&\sim&
%p_{On} \frac{r - d}{d - \rho} \pr{\DD^r > \frac{x}{r - d}}
%\hspace*{.4in} \mbox{ as } x \to \infty \\
%&=&
%\frac{\nu \delta}{1 + \nu \delta} \frac{\lambda \beta - s_L}{\bar{s} - \lambda \beta}
%\frac{C_D}{\Gamma(2 - \nu_D) \delta} \left(\frac{x}{\lambda \beta - s_L}\right)^{1 - \nu_D} \\
%&=&
%\frac{\nu \delta}{1 + \nu \delta} \frac{(\lambda \beta - s_L)^{\nu_D}}{\bar{s} - \lambda \beta}
%\frac{C_D}{\Gamma(2 - \nu_D) \delta} x^{1 - \nu_D},
%\end{eqnarray*}
%upon noting that $d - \rho = s_H - \lambda \beta - \frac{\nu \delta}{1 + \nu \delta} (s_H - s_L) = \frac{s_H + \nu \delta s_L}{1 + \nu \delta} - \lambda \beta = \bar{s} - \lambda \beta$.
Thus, we deduce that $\pr{\VV_{H, \lambda} > x} \sim \pr{\WW_{On}^d > x}$ as $x \to \infty$.
Recall again that $\VV_{H, \lambda}$ may also be interpreted as the workload in the dual-input system mentioned earlier, given that the On-Off source in that system is Off.
We infer that in case~(b) the workload asymptotics of the dual-input system are the same as in Reference System~B, which represents a reduced-load equivalence property for the workload asymptotics of the dual-input system when the gradual input is `heavier-tailed' than the instantaneous input.

Moreover, it follows from similar arguments as in~\cite{BD98a} that $\WW_{On}^d$ behaves as the workload $\VV_{\hat\lambda}^1$ in an M/G/1 queue with arrival rate $\hat\lambda = \nu / d$ and generic service requirement $\hat{\BB} = (\lambda \beta - s_L) \DD$, so that $C_{\hat{B}} = (\lambda \beta - s_L)^{\nu_D} C_D$ and $\hat\beta = (\lambda \beta - s_L) \delta$, yielding with $c = 1$,
%\begin{eqnarray*}
%\frac{\hat\lambda C_{\hat{B}}}{c - \hat\lambda \hat\beta}
%&=&
%\frac{\frac{\nu}{s_H - \lambda \beta} (\lambda \beta - s_L)^{\nu_D} C_D}{1 - \frac{\nu}{s_H - \lambda \beta} (\lambda \beta - s_L) \delta} \\
%&=&
%\frac{\nu (\lambda \beta - s_L)^{\nu_D} C_D}{s_H - \lambda \beta - \nu (\lambda \beta - s_L) \delta} \\
%&=&
%\frac{\nu (\lambda \beta - s_L)^{\nu_D} C_D}{s_H + \nu \delta s_L - (1 + \nu \delta) \lambda \beta} \\
%&=&
%\frac{\nu (\lambda \beta - s_L)^{\nu_D} C_D}{(1 + \nu \delta) (\bar{s} - \lambda \beta)} \\
%&=&
%C_{V_{H,\lambda}}
%\end{eqnarray*}
\[
\frac{\hat\lambda C_{\hat{B}}}{c - \hat\lambda \hat\beta} =
\frac{\frac{\nu}{s_H - \lambda \beta} (\lambda \beta - s_L)^{\nu_D} C_D}{1 - \frac{\nu}{s_H - \lambda \beta} (\lambda \beta - s_L) \delta} =
\frac{\nu (\lambda \beta - s_L)^{\nu_D} C_D}{(1 + \nu \delta) (\bar{s} - \lambda \beta)} =
C_{V_{H,\lambda}}
\]
as in~(\ref{cvhl1}), confirming the equivalence in the tail asymptotics.

\subsection{Heavy-traffic limit}

We now turn to the heavy-traffic limit of the workload, and assume that, for every $\omega \geq 0$,
\[
\frac{1}{\lambda^* - \lambda} \left[\Psi_\lambda(0) - \Psi_\lambda\left(\kappa_\lambda^{- \frac{1}{\alpha}} \omega\right)\right] \to 0
\hspace*{.4in} \mbox{ as } \lambda \uparrow \lambda^*,
\]
so that $G_\lambda\left(\kappa_\lambda^{- \frac{1}{\alpha}} \omega\right) \to 1$ as $\lambda \uparrow \lambda^*$.
Also, it is easily verified that $H_\lambda\left(\kappa_\lambda^{- \frac{1}{\alpha}} \omega\right) \to \omega^\alpha$ as $\lambda \uparrow \lambda^*$.
Statement~(ii) of Proposition~\ref{main1} then implies that when $\kappa_\lambda S_\lambda^\alpha \to 1$ as $\lambda \uparrow \lambda^*$,
\[
S_\lambda \VV_{H, \lambda} \to \WW
\hspace*{.4in} \mbox{ as } \lambda \uparrow \lambda^*,
\]
where $\WW$ is a random variable with a Mittag-Leffler distribution with parameter~$\alpha$.

Also, noting that
\begin{equation}
\frac{\lambda}{(\lambda^* - \lambda) \kappa_\lambda} \frac{C_{\hat{B}}}{\hat\beta} \to \zeta =
%\frac{\lambda^* - \lambda}{\lambda^*} \frac{\eta_\lambda}{(1 + \nu \delta) (\lambda - \lambda^*}) \beta} = \frac{\eta_\lambda}{(1 + \nu \delta) \lambda^* \beta} \to
\frac{(1 + \nu \delta) \bar{s}}{\eta_{\lambda^*}} \frac{C_{\hat{B}}}{\hat\beta}
\hspace*{.4in} \mbox{ as } \lambda \uparrow \lambda^*,
\label{zeta1}
\end{equation}
we find that $\VV_{H, \lambda}$ has the same heavy-traffic limit for $\lambda \uparrow \lambda^*$, up to a relative capacity slack factor~$\zeta$, as the stationary workload $\VV_{\hat\lambda}^c$ defined in Section~\ref{modelsetup} with a regularly varying service requirement distribution of index~$- \alpha$ as in~(\ref{sere1}) for $\hat\lambda \uparrow \hat\lambda^* = c / \beta$.

In case~(a), $\alpha = \nu_B - 1$ and $\eta_\lambda = \lambda (1 + \nu \delta) C_B$, so that
$\kappa_\lambda = \frac{\lambda}{\lambda^* - \lambda} \frac{C_B}{\beta}$, $\zeta = \frac{\bar{s}}{\lambda^* C_B} \frac{C_{\hat{B}}}{\hat\beta} = \frac{\beta}{C_B} \frac{C_{\hat{B}}}{\hat\beta}$, and $\kappa_\lambda \frac{\lambda^* - \lambda}{\lambda^*} \frac{\beta}{C_B} = \frac{\lambda}{\lambda^*} \to 1$ as $\lambda \uparrow \lambda^*$.
Thus,
\[
\left(\frac{\bar{s} - \rho}{\bar{s}} \frac{\beta}{C_B}\right)^{\frac{1}{\nu_B - 1}} \VV_{H, \lambda} =
%\left(\frac{(\lambda^* - \lambda) \beta}{\lambda^* C_B}\right)^{\frac{1}{\nu_B - 1}} \VV_{H, \lambda} =
%\left(\frac{\lambda}{\lambda^*} \frac{(\lambda^* - \lambda) \beta}{\lambda C_B}\right)^{\frac{1}{\nu_B - 1}} \VV_{H, \lambda} =
%\left(\frac{\lambda}{\lambda^*}\right)^{\frac{1}{\nu_B - 1}} \kappa_\lambda^{- \frac{1}{\nu_B - 1}} \VV_{H, \lambda}
\left(\frac{\lambda^* - \lambda}{\lambda^*} \frac{\beta}{C_B}\right)^{\frac{1}{\nu_B - 1}} \VV_{H, \lambda} \to \WW
\hspace*{.4in} \mbox{ as } \lambda \uparrow \lambda^*.
\]

When $\hat{\BB} = \BB$, so that $\zeta = 1$, we find that $\VV_{H, \lambda}$ has the same heavy-traffic limit for $\lambda \uparrow \lambda^* = \bar{s} / \beta$ as $\VV_{\hat\lambda}^c$ for $\hat\lambda \uparrow \hat\lambda^* = c / \beta$, as is confirmed by comparison with Example~$4$.
Specifically, if we define
\[
\lambda(\epsilon) = \lambda^* \left(1 - \frac{C_{\hat{B}}}{\hat\beta} \epsilon^{\nu - 1}\right) =
\lambda^* \left(1 - \frac{C_B}{\beta} \epsilon^{\nu - 1}\right)
\]
and
\[
\hat\lambda(\epsilon) = \hat\lambda^* \left(1 - \frac{C_{\hat{B}}}{\hat\beta} \epsilon^{\nu - 1}\right),
\]
then
\[
\lim_{\epsilon \downarrow 0} \pr{\epsilon \VV_{\lambda(\epsilon)} \leq x} =
\lim_{\epsilon \downarrow 0} \pr{\epsilon \VV_{\hat\lambda(\epsilon)}^c \leq x} = \pr{\WW \leq x}.
\]
This holds in particular for $c = \bar{s}$, in which case the workload $\VV_{\hat\lambda}^c$ corresponds to that in Reference System~A.
Further recall that $\VV_{H, \lambda}$ may be interpreted as the workload in the dual-input system, given that the On-Off source in that system is Off.
We thus conclude that in case~(a) the heavy-traffic limit of the workload in the dual-input system is the same as in Reference system~A with only gradual input, whose service speed~$\bar{s}$ equals the service rate $s_H$ of the dual-input system reduced by the time-average rate $\frac{s_H - s_L}{1 + \nu \delta}$ of the On-Off source as observed before.
The identical heavy-traffic limit in the two systems is a manifestation of a reduced-load equivalence for the dual-input system in heavy-traffic conditions when the instantaneous input is `heavier-tailed' than the gradual input.

In case~(b), $\alpha = \nu_D - 1$ and $\eta_\lambda = \nu (\lambda \beta - s_L)^{\nu_D} C_D$, so that
\[
\kappa_\lambda = \frac{\nu (\lambda \beta - s_L)^{\nu_D} C_D}{(1 + \nu \delta) (\lambda^* - \lambda) \beta},
\]
and
\[
\kappa_\lambda \frac{(\lambda^* - \lambda) (1 + \nu \delta) \beta}{\nu (\bar{s} - s_L)^{\nu_D} C_D} =
\left(\frac{\lambda \beta - s_L}{\bar{s} - s_L}\right)^{\nu_D} \to 1
\hspace*{.4in} \mbox{ as } \lambda \uparrow \lambda^*.
\]
Thus,
\[
\left(\frac{(\bar{s} - \lambda \beta) (1 + \nu \delta)}{\nu (\bar{s} - s_L)^{\nu_D} C_D}\right)^{\frac{1}{\nu_D - 1}} \VV_{H, \lambda} =
%\left(\frac{\lambda}{\lambda^*} \left(\frac{\lambda \beta - s_L}{\bar{s} - s_L}\right)^\nu \right)^{\frac{1}{\nu_D - 1}} \kappa_\lambda^{- \frac{1}{\nu_D - 1}} \VV_{H, \lambda}
%\left(\frac{(\lambda^* - \lambda) (1 + \nu \delta) \beta}{\nu C_D \lambda^* (\bar{s} - s_L)^{\nu_D}}\right)^{\frac{1}{\nu_D - 1}} \VV_{H, \lambda}
\left(\frac{(\lambda^* - \lambda) (1 + \nu \delta) \beta}{\nu (\bar{s} - s_L)^{\nu_D} C_D}\right)^{\frac{1}{\nu_D - 1}} \VV_{H, \lambda} \to \WW
\hspace*{.4in} \mbox{ as } \lambda \uparrow \lambda^*.
\]

Now let us consider Reference System~B.
As mentioned above, the workload $\WW_{On}^d$ in this system at the start of an On-period behaves as the workload $\VV_{\hat\lambda}^1$ in an M/G/1 queue with arrival rate $\hat\lambda = \nu / d$ and generic service requirement $\hat{\BB} = (\lambda \beta - s_L) \DD$, so that $C_{\hat{B}} = (\lambda \beta - s_L)^{\nu_D} C_D$ and $\hat\beta = (\lambda \beta - s_L) \delta$.

Just like in Example~$3$, we obtain
\[
\left((1 - \hat\lambda \hat\beta) \frac{\hat\beta}{C_{\hat{B}}}\right)^{\frac{1}{\nu_D - 1}} \VV_{\hat\lambda}^1 =
\left((1 - \hat\lambda (\lambda \beta - s_L) \delta) \frac{\delta}{(\lambda \beta - s_L)^{\nu_D - 1} C_D}\right)^{\frac{1}{\nu_D - 1}} \VV_{\hat\lambda}^1 \to \WW
%\left(\frac{(\nu^* - \nu) (\lambda \beta - s_L) \delta}{s_H - \lambda \beta} \frac{\delta}{(\lambda \beta - s_L)^{\nu_D - 1} C_D}\right)^{\frac{1}{\nu_D - 1}} \WW_{On}^d
\]
as $\hat\lambda \uparrow \hat\lambda^* = \frac{1}{(\lambda \beta - s_L) \delta}$.
Now observe that
\[
1 - \hat\lambda (\lambda \beta - s_L) \delta = (\hat\lambda^* - \hat\lambda) (\lambda \beta - s_L) \delta,
\]
so that
\[
(1 - \hat\lambda (\lambda \beta - s_L) \delta) \frac{\delta}{(\lambda \beta - s_L)^{\nu_D - 1} C_D} =
\frac{(\hat\lambda^* - \hat\lambda) (\lambda \beta - s_L) \delta^2}{(\bar{s} - s_L)^{\nu_D - 1} C_D}
\left(\frac{\bar{s} - s_L}{\lambda \beta - s_L}\right)^{\nu_D - 1},
%\frac{(\bar{s} - \lambda \beta) (1 + \nu \delta)}{\nu (\bar{s} - s_L)^{\nu_D} C_D} \frac{\nu \delta (\lambda \beta - s_L)}{s_H - \lambda \beta} \left(\frac{\bar{s} - s_L}{\lambda \beta - s_L}\right)^{\nu_D},
\]
with $\frac{\bar{s} - s_L}{\lambda \beta - s_L} \to 1$ as $\lambda \uparrow \lambda^*$.
Thus, we deduce that
\[
%\left(\frac{(\bar{s} - \lambda \beta) (1 + \nu \delta)}{\nu (\bar{s} - s_L)^{\nu_D} C_D}\right)^{\frac{1}{\nu_D - 1}} \WW_{On}^d
\left(\frac{(\hat\lambda^* - \hat\lambda) (\lambda \beta - s_L) \delta^2}{(\bar{s} - s_L)^{\nu_D - 1} C_D}\right)^{\frac{1}{\nu_D - 1}} \VV_{\hat\lambda}^1 \to \WW
\hspace*{.4in} \mbox{ as } \hat\lambda \uparrow \hat\lambda^*,
\]
indicating that $\VV_{H,\lambda}$ has the same heavy-traffic limit for $\lambda \uparrow \lambda^*$, up to a relative capacity slack factor
\[
\zeta = \frac{(1 + \nu \delta) \lambda^* \beta}{\nu (\bar{s} - s_L)^{\nu_D} C_D} \left[\frac{(\lambda \beta - s_L) \hat\lambda^* \delta}{(\bar{s} - s_L)^{\nu_D - 1} C_D}\right]^{- 1} = \frac{(1 + \nu \delta) \bar{s}}{\nu \delta (\bar{s} - s_L)},
\]
as $\WW_{On}^d = \VV_{\hat\lambda}^1$ for $\hat\lambda \uparrow \hat\lambda^*$, where the dependence of $\WW_{On}^d$ on~$\hat\lambda$ is not explicitly reflected in the notation.
Since $\eta_{\lambda^*} = \nu (\lambda^* \beta - s_L)^{\nu_D} C_D = \nu (\bar{s} - s_L)^{\nu_D} C_D$ the above expression for~$\zeta$ agrees with~(\ref{zeta1}) confirming the equivalence of the heavy-traffic limit.
Specifically, if we define
\[
\lambda(\epsilon) =
\lambda^* \left(1 - \frac{1}{\zeta} \frac{(\bar{s} - s_L)^{\nu_D - 1} C_D}{(\lambda \beta - s_L) \delta \hat\lambda^*} \epsilon^{\nu_D - 1}\right) =
\lambda^* \left(1 - \frac{\nu \delta (\bar{s} - s_L)^{\nu_D} C_D}{(1 + \nu \delta) \bar{s}} \epsilon^{\nu_D - 1}\right)
%\lambda^* \left(1 - \frac{\nu \delta (\bar{s} - s_L)^{\nu_D} C_D}{(1 + \nu \delta) \lambda^* \beta} \epsilon^{\nu_D - 1}\right) =
%\lambda^* - \frac{\nu \delta (\bar{s} - s_L)^{\nu_D} C_D}{(1 + \nu \delta) \beta} \epsilon^{\nu_D - 1}
\]
and
\[
\hat\lambda(\epsilon) = \hat\lambda^* \left(1 - \frac{(\bar{s} - s_L)^{\nu_D - 1} C_D}{(\lambda \beta - s_L) \delta \hat\lambda^*} \epsilon^{\nu_D - 1}\right) =
\hat\lambda^* \left(1 - (\bar{s} - s_L)^{\nu_D - 1} C_D \epsilon^{\nu_D - 1}\right),
\]
then
\[
\lim_{\epsilon \downarrow 0} \pr{\epsilon \VV_{\lambda(\epsilon)} \leq x} =
\lim_{\epsilon \downarrow 0} \pr{\epsilon \VV_{\hat\lambda(\epsilon)}^1 \leq x} = \pr{\WW \leq x}.
\]

Recall again that $\VV_{H, \lambda}$ may also be interpreted as the workload in the dual-input system, given that the On-Off source in that system is Off.
We thus conclude that in case~(b) the heavy-traffic limit of the workload in the dual-input system is the same as in Reference System~B, whose service speed $d = s_H - \lambda \beta$ equals the service rate of the dual-input system reduced by the average load $\lambda \beta$ generated by the customers.
The identical heavy-traffic limits in the two systems represents a reduced-load equivalence property for the dual-input system when the gradual input is `heavier-tailed' than the instantaneous input.

\section{M/G/2 queue with heterogeneous servers}
\label{heterogeneous}

In this section we consider an M/G/2 queue with heterogeneous servers as studied by Boxma, Deng \& Zwart~\cite{BDZ02}.
Customers arrive as a Poisson process of rate~$\lambda$ and are served in a FCFS manner.
The service times are exponentially distributed with parameter~$\mu$ at server~$1$, and are i.i.d.\ copies of a non-negative random variable~$\BB$ with LST $\beta(\omega)$ and mean $\beta < \infty$ at server~$2$.
It is easily seen that the stability condition is $\lambda < \lambda^* = \mu + \frac{1}{\beta}$, which is henceforth assumed to be satisfied.
In addition, we will assume that $\lambda > \mu$, implying that server~$1$ alone cannot handle the entire workload and that server~$2$ is thus needed to achieve stability. \\

Let $\WW_\lambda$ be the waiting time of an arbitrary customer in equilibrium.
Denote by $\pi_0$ the probability that the system is empty in equilibrium so that both servers are idle, by $\pi_1$ the probability that there is one customer in the system being served by server~$1$ so that server~$2$ is idle, and by $\pi_2$ the probability that there is one customer in the system being served by server~$2$ so that server~$1$ is idle.
Then the time-average service rate may be expressed as
\[
(1 - \pi_0 - \pi_2) \mu + (1 - \pi_0 - \pi_1) \frac{1}{\beta},
\]
which must be equal to the arrival rate~$\lambda$ in case the system is stable, yielding the identity relation
\begin{equation}
\frac{1}{\beta} \pi_0 + \mu \pi_0 + \frac{1}{\beta} \pi_1 + \mu \pi_2 = \frac{1}{\beta} + \mu - \lambda
\label{identity1}
\end{equation}
as stated in Equation (2.14) in~\cite{BDZ02}.

Also, it is shown in~\cite{BDZ02} that
\[
\expect{\ee^{- \omega \WW_\lambda}} = F_\lambda(\omega) + \frac{P_\lambda(\omega)}{Q_\lambda(\omega)},
\]
with
\begin{eqnarray*}
F_\lambda(\omega)
&=&
\pi_0 + \pi_1 - \frac{\mu + \omega}{\lambda - \mu - \omega} \pi_2, \\
P_\lambda(\omega)
&=&
\hat{h}\left(\frac{\mu + \hat{g}(\omega)}{\lambda}\right) - \frac{\mu + \hat{g}(\omega)}{\lambda} \left[\lambda \pi_0 + \hat{g}(\omega) \pi_1\right] \\
&+&
\left[(\lambda - \mu - \omega) \left[\lambda \pi_0 + (\lambda - \mu - \omega) \pi_1\right] \beta -
\frac{\lambda \beta \hat{g}(\omega)}{\mu + \hat{g}(\omega)}
\hat{h}\left(\frac{\mu + \hat{g}(\omega)}{\lambda}\right)\right] \beta^r(\hat{f}(\omega)), \\
Q_\lambda(\omega)
&=&
(\lambda - \mu - \omega) \left[1 - \frac{\omega}{\lambda} - (\lambda - \mu - \omega) \beta \beta^r(\hat{f}(\omega))\right], \\
\hat{f}(\omega)
&=&
\frac{\omega (\lambda - \mu - \omega)}{\lambda - \omega}, \\
\hat{g}(\omega)
&=&
\omega - \hat{f}(\omega) = \frac{\mu \omega}{\lambda - \omega},
\end{eqnarray*}
and $\hat{h}(\cdot)$ as defined in~\cite{BDZ02}.

Note that
\begin{equation}
F_\lambda(\omega) - F_\lambda(0) =
\left(\frac{\mu}{\lambda - \mu} - \frac{\mu + \omega}{\lambda - \mu - \omega}\right) \pi_2 \leq f \omega
\label{flomega1}
\end{equation}
for some finite constant~$f$ independent of~$\lambda$, and that the expression for $\expect{\ee^{- \omega \WW_\lambda}}$ matches the form of~(\ref{form1}), with
\[
G_\lambda(\omega) = \frac{P_\lambda(\omega)}{(\lambda - \mu) [1 - (\lambda - \mu) \beta]},
\]
and
\[
H_\lambda(\omega) = \frac{Q_\lambda(\omega)}{(\lambda - \mu) [1 - (\lambda - \mu) \beta]} - 1.
\]

Straightforward algebraic manipulations yield
\[
\hspace*{-.2in}
H_\lambda(\omega)=
\frac{1}{\lambda^* - \lambda} \frac{1}{\lambda - \mu} \left[\left[\frac{\mu}{\lambda \beta} - 2 [\frac{1}{\beta} - (\lambda - \mu)] + \omega (\frac{1}{\lambda \beta} - 1)\right] \omega + (\lambda - \mu - \omega)^2 (1 - \beta^r(\hat{f}(\omega)))\right].
\]

It follows from the results in~\cite{BDZ02} that if the service requirement distribution at server~$2$ is regularly varying of index $\nu \in (1, 2)$ as in~(\ref{sere2}), then
\begin{eqnarray*}
& &
(\lambda - \mu) [1 - (\lambda - \mu) \beta] G_\lambda(\omega) - (\lambda - \mu) (\lambda \pi_0 + \lambda \pi_1 - \mu \pi_1) \beta + \mu \pi_2 \\
&=&
P_\lambda(\omega) - (\lambda - \mu) (\lambda \pi_0 + \lambda \pi_1 - \mu \pi_1) \beta + \mu \pi_2 \\
&=&
- (\lambda - \mu) (\lambda \pi_0 + \lambda \pi_1 - \mu \pi_1) C_B \left(\frac{\lambda - \mu}{\lambda}\right)^{\nu - 1} \omega^{\nu - 1} + \oo(\omega^{\nu - 1})
\hspace*{.4in} \mbox{ as } \omega \downarrow 0.
\end{eqnarray*}
Thus,
\begin{equation}
G_\lambda(0) - G_\lambda(\infty) = \gamma_\lambda \omega^\alpha + \oo(\omega^\alpha)
\hspace*{.4in} \mbox{ as } \omega \downarrow 0,
\label{glomega1}
\end{equation}
with $\alpha = \nu - 1 \in (0, 1)$,
\begin{equation}
G_\lambda(0) =
\frac{(\lambda - \mu) (\lambda \pi_0 + \lambda \pi_1 - \mu \pi_1) \beta + \mu \pi_2}{(\lambda - \mu) [1 - (\lambda - \mu) \beta]},
\label{glzero1}
\end{equation}
and
\begin{equation}
\gamma_\lambda =
%\frac{(\lambda - \mu) (\lambda Q_0 + \lambda Q_1 - \mu Q_1) C_B}{(\lambda - \mu) [1 - (\lambda - \mu) \beta]} \left(\frac{\lambda - \mu}{\lambda}\right)^{\nu - 1}
\frac{(\lambda \pi_0 + \lambda \pi_1 - \mu \pi_1) C_B}{1 - (\lambda - \mu) \beta} \left(\frac{\lambda - \mu}{\lambda}\right)^{\nu - 1}.
\label{gammal1}
\end{equation}

Also,
\begin{eqnarray*}
& &
(\lambda - \mu) [1 - (\lambda - \mu) \beta] H_\lambda(\omega) \\
&=&
(\lambda - \mu) [1 - (\lambda - \mu) \beta] [H_\lambda(\omega) + 1] - (\lambda - \mu) [1 - (\lambda - \mu) \beta] \\
&=&
Q_\lambda(\omega) - (\lambda - \mu) [1 - (\lambda - \mu) \beta] \\
&=&
\left((\lambda - \mu) [1 - (\lambda - \mu) \beta]\right)^2 \frac{C_B}{[1 - (\lambda - \mu) \beta]^2} \left(\frac{\lambda - \mu}{\lambda}\right)^{\nu - 1} \omega^{\nu - 1} + \oo(\omega^{\nu - 1}) \\
&=&
(\lambda - \mu)^2 C_B \left(\frac{\lambda - \mu}{\lambda}\right)^{\nu - 1} \omega^{\nu - 1} + \oo(\omega^{\nu - 1})
\hspace*{.4in} \mbox{ as } \omega \downarrow 0.
\end{eqnarray*}

Thus,
\begin{equation}
H_\lambda(\omega) = \kappa_\lambda \omega^\alpha + \oo(\omega^\alpha)
\hspace*{.4in} \mbox{ as } \omega \downarrow 0,
\label{hlomega1}
\end{equation}
with $\alpha = \nu - 1 \in (0, 1)$ and
\begin{equation}
\kappa_\lambda =
%\frac{(\lambda - \mu)^2 C_B}{(\lambda - \mu) [1 - (\lambda - \mu) \beta]}  \left(\frac{\lambda - \mu}{\lambda}\right)^{\nu - 1}
\frac{(\lambda - \mu) C_B}{1 - (\lambda - \mu) \beta} \left(\frac{\lambda - \mu}{\lambda}\right)^{\nu - 1}.
\label{kappal1}
\end{equation}

\iffalse

\[
\hat{H}^\alpha(\lambda, \omega) \to \hat\kappa_\lambda
\hspace*{.4in} \mbox{ as } \omega \downarrow 0
\]
for all $\lambda \in [0, \lambda^*]$, and
\[
\hat{H}\left(\lambda, \left(\frac{\lambda^* - \lambda}{\gamma \lambda^*}\right)^{\frac{1}{\alpha}}\right) \to \hat\kappa_{\lambda^*}
\hspace*{.4in} \mbox{ as } \lambda \uparrow \lambda^*
\]
for every $\omega \geq 0$, with
\begin{eqnarray*}
\hat\kappa_\lambda
&=&
\frac{\lambda^* - \lambda}{\lambda} \kappa_\lambda \\
&=&
\frac{\lambda^* - \lambda}{\lambda} \frac{(\lambda - \mu) C_B}{1 - (\lambda - \mu) \beta} \left(\frac{\lambda - \mu}{\lambda}\right)^{\nu - 1} \\
&=&
\frac{\lambda^* - \lambda}{\lambda} \frac{(\lambda - \mu) \frac{C_B}{\beta}}{\frac{1}{\beta} - \lambda + \mu} \left(\frac{\lambda - \mu}{\lambda}\right)^{\nu - 1} \\
&=&
\frac{\lambda^* - \lambda}{\lambda}
\frac{(\lambda - \mu) C_B}{(\lambda^* - \lambda) \beta} \left(\frac{\lambda - \mu}{\lambda}\right)^{\nu - 1} \\
&=&
\frac{(\lambda - \mu) C_B}{\lambda \beta} \left(\frac{\lambda - \mu}{\lambda}\right)^{\nu - 1} \\
&=&
\frac{C_B}{\beta} \left(\frac{\lambda - \mu}{\lambda}\right)^\nu,
\end{eqnarray*}

\[
\gamma = \hat\kappa_{\lambda^*} = \frac{(\lambda^* - \mu) C_B}{\lambda^* \beta} \left(\frac{\lambda^* - \mu}{\lambda^*}\right)^{\nu - 1} = \frac{C_B}{\beta} \left(\frac{\lambda^* - \mu}{\lambda^*}\right)^\nu.
\]

\fi

\subsection{Tail asymptotics}

We first consider the tail asymptotics of the waiting time as obtained in~\cite{BDZ02}, but we approach these from the perspective of~(\ref{form1}) in order to highlight the connection with the heavy-traffic limit.
Invoking (\ref{flomega1},\ref{glomega1},\ref{glzero1},\ref{gammal1},\ref{hlomega1},\ref{kappal1}), and noting that $F_\lambda(\omega) - F_\lambda(0) = \oo(\omega^{\nu - 1})$, statement~(i) of Proposition~\ref{main1} implies that
\begin{eqnarray}
\pr{\WW_\lambda > x} \sim \frac{C_{W_\lambda}}{\Gamma(2 - \nu)} x^{1 - \nu}
&=&
\frac{1 - \pi_0 - \pi_1}{1 - (\lambda - \mu) \beta} \frac{C_B}{\Gamma(2 - \nu) \beta} \left(\frac{\lambda x}{\lambda - \mu}\right)^{1 - \nu} \nonumber \\
&\sim&
\frac{1 - \pi_0 - \pi_1}{1 - (\lambda - \mu) \beta} \pr{\BB^r > \frac{\lambda x}{\lambda - \mu}}
\hspace*{.4in} \mbox{ as } x \to \infty,
\label{tailwlx}
\end{eqnarray}
with
{\small
\begin{eqnarray*}
\!\!\!\!\!\!\!\!\!\!\!\!\!\!\!\!\!\!\!\!\!\!\!\!\!\!\!\!\!\!\!\!\!\!\!\!\!\!\!\!\!\!\!\!\!\!\!\!\!\!\!\!\!\!
\!\!\!\!\!\!\!\!\!\!\!\!\!\!\!\!\!\!\!\!\!\!\!\!\!\!\!\!\!\!\!\!\!\!\!\!\!\!\!\!\!\!\!\!\!\!\!\!\!\!\!\!\!\!
C_{W_\lambda}
&=&
\gamma_\lambda + \kappa_\lambda G_\lambda(0) \\
&=&
\frac{(\lambda \pi_0 + \lambda \pi_1 - \mu \pi_1) C_B}{1 - (\lambda - \mu) \beta} \left(\frac{\lambda - \mu}{\lambda}\right)^{\nu - 1} +
\frac{(\lambda - \mu) (\lambda \pi_0 + \lambda \pi_1 - \mu \pi_1) \beta + \mu \pi_2}{(\lambda - \mu) [1 - (\lambda - \mu) \beta]} \frac{(\lambda - \mu) C_B}{1 - (\lambda - \mu) \beta} \left(\frac{\lambda - \mu}{\lambda}\right)^{\nu - 1} \\
&=&
\left[\frac{[1 - (\lambda - \mu) \beta] (\lambda \pi_0 + \lambda \pi_1 - \mu \pi_1)}{1 - (\lambda - \mu) \beta} + \frac{(\lambda - \mu) (\lambda \pi_0 + \lambda \pi_1 - \mu \pi_1) \beta + \mu \pi_2}{1 - (\lambda - \mu) \beta}\right]
\frac{C_B}{1 - (\lambda - \mu) \beta} \left(\frac{\lambda - \mu}{\lambda}\right)^{\nu - 1} \\
&=&
\frac{\lambda \pi_0 + \lambda \pi_1 - \mu \pi_1 + \mu \pi_2}{[1 - (\lambda - \mu) \beta]^2}
\frac{C_B}{\beta} \left(\frac{\lambda - \mu}{\lambda}\right)^{\nu - 1} \\
&=&
\frac{1 - \pi_0 - \pi_1}{1 - (\lambda - \mu) \beta} \frac{C_B}{\beta} \left(\frac{\lambda - \mu}{\lambda}\right)^{\nu - 1},
\end{eqnarray*}
}
where the last step follows from the identity relation~(\ref{identity1}).
This is consistent with Theorem~4.1 in~\cite{BDZ02}.

Also, $\WW_\lambda$ has the same tail asymptotics as the stationary workload $\VV_{\hat\lambda}^c$ defined in Section~\ref{modelsetup} with a regularly varying service requirement distribution of index~$- \nu$ as in~(\ref{sere1}) when $C_{W_\lambda} = \frac{\hat\lambda C_{\hat{B}}}{c - \hat\lambda \hat\beta} = \frac{\hat\lambda \hat\beta}{c - \hat\lambda \hat\beta} \frac{C_{\hat{B}}}{\hat\beta}$.
In case $\hat{\BB} = \frac{\lambda}{\lambda - \mu} \BB$, so that $C_{\hat{B}} = \left(\frac{\lambda}{\lambda - \mu}\right)^{- \nu}$ and $\hat\beta = \frac{\lambda}{\lambda - \mu} \beta$, this equation reduces to
\begin{equation}
\frac{1 - \pi_0 - \pi_1}{1 - (\lambda - \mu) \beta} = \frac{\hat\lambda \hat\beta}{c - \hat\lambda \hat\beta}.
\label{redu1}
\end{equation}
This may be intuitively interpreted as follows.
First of all, it is important to observe that it is \textit{not} the case that the waiting time in the system with two heterogeneous servers roughly behaves as the workload in a single-server queue.
In particular, while the waiting time gradually grows large at a rate $\lambda / \mu - 1 > 0$ when a large service time occurs at server~$2$, the workload in a single-server queue instantaneously jumps to a large value when a customer with a large service requirement arrives.
The explanation why $\WW_\lambda$ and $\VV_{\hat\lambda}^c$ nevertheless have the same tail asymptotics, is because the waiting time in the system with two heterogeneous servers behaves approximately as the workload $\WW^d$ in a fluid queue which increases at a rate $\lambda / \mu - 1$ during On-periods, with the On-periods having the distribution of the service times at server~$2$, scaled by a factor $\mu / \lambda < 1$.
The average number of On-periods per time unit equals the average number of service completions at server~$2$ per time unit, which is $(1 - \pi_0 - \pi_1) / \beta$, and the workload, when positive, decreases on average at rate
\[
1 - \lambda \frac{\frac{\beta}{\mu}}{\beta + \frac{1}{\mu}}.
\]

In terms of the notation introduced in Section~\ref{modelsetup}, we have a fluid queue with generic On-period $\AAAA = \frac{\mu}{\lambda} \BB$, $1 / (\expect{\AAAA} + \expect{\UU}) = (1 - \pi_0 - \pi_1) / \beta$, $r - d = \frac{\lambda}{\mu} - 1$, and $d - \rho = 1 - \lambda \frac{\beta}{1 + \beta \mu}$, so that
\[
r - \rho = \frac{\lambda}{\mu} - \lambda \frac{\frac{\beta}{\mu}}{\beta + \frac{1}{\mu}} = \frac{\frac{\lambda}{\mu}}{1 + \beta \mu},
\]
and
\[
\frac{r - \rho}{d - \rho} = \frac{\frac{\frac{\lambda}{\mu}}{1 + \beta \mu}}{1 - \lambda \frac{\beta}{1 + \beta \mu}} =
\frac{\frac{\lambda}{\mu}}{\beta \mu + 1 - \lambda \beta} =
\frac{\frac{\lambda}{\mu}}{1 - (\lambda - \mu) \beta}.
\]
In view of~(\ref{tailwcx}), this heuristic observation yields
\begin{eqnarray*}
\pr{\WW_\lambda > x}
&\sim&
\pr{\WW^d > x} \\
&\sim&
%p_{On} \frac{r - \rho}{d - \rho} \pr{\AAAA^r > \frac{x}{r - d}} \\
%&=&
\frac{\expect{\AAAA}}{\expect{\AAAA} + \expect{\UU}} \frac{r - \rho}{d - \rho} \pr{\AAAA^r > \frac{x}{r - d}} \\
&=&
\frac{\mu}{\lambda} \beta \frac{1 - \pi_0 - \pi_1}{\beta} \frac{\frac{\lambda}{\mu}}{1 - (\lambda - \mu) \beta} \pr{\frac{\mu}{\lambda} \BB^r > \frac{x}{\frac{\lambda}{\mu} - 1}} \\
%&=&
%\frac{1 - \pi_0 - \pi_1}{1 - (\lambda - \mu) \beta} \pr{\BB^r > \frac{\lambda x}{\lambda - \mu}} \\
&\sim&
\frac{1 - \pi_0 - \pi_1}{1 - (\lambda - \mu) \beta} \frac{C_B}{\Gamma(2 - \nu)} \left(\frac{\lambda x}{\lambda - \mu}\right)^{1 - \nu} \\
&=&
\frac{C_{W_\lambda}}{\Gamma(2 - \nu)} x^{1 - \nu}
\hspace*{.4in} \mbox{ as } x \to \infty,
\end{eqnarray*}
which is in agreement with~(\ref{tailwlx}).
As mentioned in Section~\ref{modelsetup}, when $\AAAA = \frac{\mu}{\lambda} \BB$ and $\hat{\BB} = \frac{\lambda - \mu}{\lambda} \BB$ so that $\hat{\BB} = \frac{\lambda - \mu}{\mu} \AAAA = (r - d) \AAAA$, the workload $\WW^d$ in such a fluid queue has the same tail asymptotics as $\VV_{\hat\lambda}^c$ when
\[
(1 - p_{On}) \frac{\rho}{d - \rho} = \frac{\hat\lambda \hat\beta}{c -\hat\lambda \hat\beta},
\]
which corresponds to~(\ref{redu1}) as derived above.

\subsection{Heavy-traffic limit}

We now turn to the heavy-traffic limit of the waiting time.
In view of~(\ref{flomega1},\ref{glomega1},\ref{kappal1}), we have
\[
F_\lambda\left(\kappa_\lambda^{- \frac{1}{\nu - 1}} \omega\right) \to 0
\hspace*{.4in} \mbox{ as } \lambda \uparrow \lambda^*
\]
and
\[
G_\lambda\left(\kappa_\lambda^{- \frac{1}{\nu - 1}} \omega\right) \to 1
\hspace*{.4in} \mbox{ as } \lambda \uparrow \lambda^*.
\]

Also, it is easily verified that
\[
H_\lambda\left(\kappa_\lambda^{- \frac{1}{\nu - 1}} \omega\right) \to \omega^{\nu - 1}
\hspace*{.4in} \mbox{ as } \lambda \uparrow \lambda^*.
\]

Statement~(ii) of Proposition~\ref{main1} then implies
\[
\kappa_\lambda^{- \frac{1}{\nu - 1}} \WW_\lambda \to \WW
\hspace*{.4in} \mbox{ as } \lambda \uparrow \lambda^*,
\]
where $\WW$ is a random variable with a Mittag-Leffler distribution with parameter $\nu - 1$.

Also, noting that
\begin{eqnarray*}
\frac{\lambda}{(\lambda^* - \lambda) \kappa_\lambda} \frac{C_{\hat{B}}}{\hat\beta}
&=&
\frac{\lambda}{\lambda^* - \lambda} \frac{1 - (\lambda - \mu) \beta}{(\lambda - \mu) C_B} \left(\frac{\lambda - \mu}{\lambda}\right)^{1 - \nu} \frac{C_{\hat{B}}}{\hat\beta} \\
&=&
\frac{\lambda}{\lambda^* - \lambda} \frac{\left(\frac{1}{\beta} + \mu - \lambda\right) \beta}{(\lambda - \mu) C_B} \left(\frac{\lambda - \mu}{\lambda}\right)^{1 - \nu} \frac{C_{\hat{B}}}{\hat\beta} \\
&=&
\frac{\lambda}{\lambda^* - \lambda} \frac{\lambda^* - \lambda}{\lambda - \mu} \frac{\beta}{C_B} \left(\frac{\lambda - \mu}{\lambda}\right)^{1 - \nu} \frac{C_{\hat{B}}}{\hat\beta} \\
&=&
\left(\frac{\lambda - \mu}{\lambda}\right)^{- \nu} \frac{\beta}{C_B} \frac{C_{\hat{B}}}{\hat\beta} \\
&=&
\left(\frac{\lambda}{\lambda - \mu}\right)^\nu \frac{\beta}{C_B} \frac{C_{\hat{B}}}{\hat\beta} \\
&\to&
\zeta \hspace*{.4in} \mbox{ as } \lambda \uparrow \lambda^*,
\end{eqnarray*}
with
\[
\zeta =
\left(\frac{\lambda^*}{\lambda^* - \mu}\right)^\nu \frac{\beta}{C_B} \frac{C_{\hat{B}}}{\hat\beta} =
\left(\frac{\mu + \frac{1}{\beta}}{\frac{1}{\beta}}\right)^\nu \frac{\beta}{C_B} \frac{C_{\hat{B}}}{\hat\beta} =
(1 + \mu \beta)^\nu \frac{\beta}{C_B} \frac{C_{\hat{B}}}{\hat\beta},
\]
we find that $\WW_\lambda$ has the same heavy-traffic limit for $\lambda \uparrow \lambda^*$, up to a relative capacity slack factor~$\zeta$, as the stationary workload $\VV_{\hat\lambda}^c$ as defined in Section~\ref{modelsetup} with a regularly varying service requirement distribution of index~$- \nu$ as in~(\ref{sere1}) for $\hat\lambda \uparrow \hat\lambda^* = c / \hat\beta$.
In case $\hat{\BB} = \frac{\lambda^* - \mu}{\lambda^*} \BB = \frac{\frac{1}{\beta}}{\mu + \frac{1}{\beta}} \BB = \frac{1}{1 + \mu \beta} \BB$, so that $C_{\hat{B}} = (1 + \mu \beta)^{- \nu} C_B$ and $\hat\beta = \beta / (1 + \mu \beta)$, we have $\zeta = 1 + \mu \beta$.
This may be heuristically explained as follows.
As before, it is important to note here that it is \textit{not} the case that the waiting time in the system with two heterogeneous servers roughly behaves as the workload in a single-server queue.
The reason why $\WW_\lambda$ and $\VV_{\hat\lambda}^c$ nevertheless have the same heavy-traffic limit, is again because the waiting time in the system with two heterogeneous servers behaves approximately as the workload $\WW^d$ in the fluid queue specified in the previous subsection.

In heavy-traffic conditions, the latter workload at an arbitrary epoch is approximately similar to the workload at the start of an On-period, which in turn evolves as the waiting time in an M/G/1 queue with arrival rate $1 / \expect{\UU}$ and generic service requirement $\hat\BB = (r - d^*) \AAAA = (\frac{\lambda^*}{\mu} - 1) \frac{\mu}{\lambda^*} \BB = \frac{\lambda^* - \mu}{\lambda^*} \BB = \frac{1}{1 + \mu \beta} \BB$.
Specifically, if we take $c = 1$ and define
\[
\lambda(\epsilon) =
\lambda^* \left(1 - \frac{1}{\zeta} \frac{\hat\beta}{C_{\hat{B}}} (1 + \mu \beta)^{- \nu} \epsilon^{\nu - 1}\right) =
%\lambda^* - \frac{\lambda^*}{\zeta} \frac{\hat\beta}{C_{\hat{B}}} (1 + \mu \beta)^{- \nu} \epsilon^{\nu - 1} =
\lambda^* - \frac{1}{\beta} \frac{\hat\beta}{C_{\hat{B}}} (1 + \mu \beta)^{- \nu} \epsilon^{\nu - 1} =
\lambda^* - \frac{1}{(1 + \mu \beta) C_B} \epsilon^{\nu - 1}
\]
and
\[
\hat\lambda(\epsilon) =
\hat\lambda^* \left(1 - \frac{\hat\beta}{C_{\hat{B}}} (1 + \mu \beta)^{- \nu} \epsilon^{\nu - 1}\right) =
%\hat\lambda^* - \hat\lambda^* \frac{\hat\beta}{C_{\hat{B}}} (1 + \mu \beta)^{- \nu} \epsilon^{\nu - 1} =
\hat\lambda^* - \frac{1}{\beta} \frac{\hat\beta}{C_{\hat{B}}} (1 + \mu \beta)^{- \nu} \epsilon^{\nu - 1} =
\hat\lambda^* - \frac{1}{(1 + \mu \beta) C_B} \epsilon^{\nu - 1},
\]
then
\[
\lim_{\epsilon \downarrow 0} \pr{\epsilon \WW_{\lambda(\epsilon)} \leq x} =
\lim_{\epsilon \downarrow 0} \pr{\epsilon \VV_{\hat\lambda(\epsilon)}^1 \leq x} = \pr{\WW \leq x}.
\]
As before, $\frac{1}{\zeta} = \frac{1}{1 + \mu \beta}$ represents the fraction of heavy-tailed load among the total load in heavy-traffic conditions, and $\lambda(\epsilon)$ and $\hat\lambda(\epsilon)$ correspond to the same absolute slack with respect to the critical arrival rates $\lambda^* = \mu + \frac{1}{\beta}$ and $\hat\lambda^* = \frac{1}{\beta}$ in the two respective systems.


\begin{thebibliography}{99}

\bibitem{AMN99}
R. Agrawal, A.M. Makowski, Ph. Nain (1999).
On a reduced-load equivalence for fluid queues under subexponentiality.
\textit{Queueing Systems} \textbf{33 (1--3)},
5--41.

\bibitem{BBRZ19}
M. Bazhba, J. Blanchet, C.-H. Rhee, B. Zwart (2019).
Queue length asymptotics for the multiple-server queue with heavy-tailed Weibull service times.
\textit{Queueing Systems} \textbf{93 (3--4)},
195--226.

\bibitem{BBRZ20}
M. Bazhba, J. Blanchet, C.-H. Rhee, B. Zwart (2020).
Sample path large deviations for L\'evy processes and random walks with Weibull increments.
\textit{Ann.\ Appl.\ Prob.} \textbf{30 (6)},
2695--2739.

%\bibitem{BRZ22}
%M. Bazhba, C.-H. Rhee, B. Zwart (2022).
%Large deviations for stochastic fluid networks with Weibullian tails.
%\textit{Queueing Systems} \textbf{102 (1--2)},
%25--52.

\bibitem{BSTW95}
J. Beran, R. Sherman, R., M.S. Taqqu, W. Willinger (1995).
Long-range dependence in variable-bit-rate video traffic.
\textit{IEEE Trans.\ Commun.} \textbf{43 (2--4)},
1566--1579.

\bibitem{BGT87}
N.H. Bingham, C.M. Goldie, J.L. Teugels (1987).
\textit{Regular Variation.}
Cambridge University Press.

\bibitem{Borovkov76}
A.A. Borovkov (1976).
\textit{Stochastic Processes in Queueing Theory}.
Springer, New York.

\bibitem{Boxma96}
O.J. Boxma (1996).
Fluid queues and regular variation.
\textit{Perf.\ Eval.} \textbf{27 \& 28},
699--712.

\bibitem{Boxma97}
O.J. Boxma (1997).
Regular variation in a multi-source fluid queue.
In: \textit{Teletraffic Contributions for the Information Age},
Proc. ITC-15. eds.\ V. Ramaswami, P. Wirth (North-Holland, Amsterdam), 391–402.

\bibitem{BC98}
O.J. Boxma, J.W. Cohen (1998).
The M/G/1 queue with heavy-tailed service time distribution.
\textit{IEEE J. Sel.\ Areas Commun.} \textbf{16},
749--763.

\bibitem{BC99a}
O.J. Boxma, J.W. Cohen (1999).
Heavy-traffic analysis for the GI/G/1 queue with heavy-tailed distributions.
\textit{Queueing Systems} \textbf{33 (1--3)},
177--204.

\bibitem{BC99b}
O.J. Boxma, J.W. Cohen (1999).
The single server queue: heavy tails and heavy traffic.
In: \textit{Self-Similar Network Traffic and Performance Evaluation},
eds.\ K. Park, W. Willinger (Wiley, New York).

\bibitem{BDZ02}
O.J. Boxma, Q. Deng, B. Zwart (2002).
Waiting-time asymptotics for the M/G/2 queue with heterogeneous servers.
\textit{Queueing Systems} \textbf{40},
5--31.

\bibitem{BD98a}
O.J. Boxma, V. Dumas (1998).
Fluid queues with long-tailed activity period distributions.
\textit{Comput.\ Commun.} \textbf{21 (17)},
1509--1529.

\bibitem{BD98b}
O.J. Boxma, V. Dumas (1998).
The busy period in the fluid queue.
\textit{Perf.\ Eval.\ Rev.} \textbf{26},
100--110.

\bibitem{BK01}
O.J. Boxma, I.A. Kurkova (2001).
The M/G/1 queue with two service rates.
\textit{Adv.\ Appl.\ Prob.} \textbf{33},
520--540.

\bibitem{CBRZ19}
B. Chen, J. Blanchet, C.-H. Rhee, B. Zwart (2019).
Efficient rare-event simulation for multiple jump events in regularly varying random walks and compound Poisson processes.
\textit{Math.\ Oper.\ Res.} \textbf{44 (3)},
919--942.

\bibitem{Cohen73}
J.W. Cohen (1973).
Some results on regular variation for distributions in queueing and fluctuation theory.
\textit{J. Appl.\ Prob.} \textbf{10},
343--353.

%\bibitem{Cohen94}
%J.W. Cohen (1994).
%On the effective bandwidth in buffer design for the multi-server channels.
%CWI Report BS-R9406.

\bibitem{Cohen97a}
J.W. Cohen (1997).
A heavy-traffic theorem for the GI/G/1 queue with a Pareto-type service time distribution.
\textit{J. Appl.\ Math.\ Stochastic Anal.} \textbf{11},
247--254.

\bibitem{Cohen97b}
J.W. Cohen (1997).
On the M/G/1 queue with heavy-tailed service time distributions.
Technical report PNA-R9702, CWI, Amsterdam.

\bibitem{Cohen97c}
J.W. Cohen (1997).
Heavy-traffic limit theorems for the heavy-tailed GI/G/1 queue.
Technical report PNA-R9719, CWI, Amsterdam.

\bibitem{Cohen98a}
J.W. Cohen (1998).
Heavy-traffic theory for the heavy-tailed M/G/1 queue and $\nu$-stable L\'evy noise traffic.
Technical report PNA-R9805, CWI, Amsterdam.

\bibitem{Cohen98b}
J.W. Cohen (1998).
The $\nu$-stable L\'evy motion in heavy-traffic analysis of queueing models with heavy-tailed distributions.
Technical report PNA-R9808, CWI, Amsterdam.

%\bibitem{CB83}
%J.W. Cohen, O.J. Boxma (1983).
%\textit{Boundary Value Problems in Queueing System Analysis.}
%North-Holland.

\bibitem{CB96}
M. Crovella, A. Bestavros (1996).
Self-similarity in World Wide Web traffic: evidence and possible causes.
In: \textit{Proc.\ ACM SIGMETRICS ’96},
160--169.

\bibitem{EV82}
P. Embrechts, N. Veraverbeke (1982).
Estimates for the probability of ruin with special emphasis on the possibility of large claims.
\textit{Insurance Math. Econom.} \textbf{1},
55--72.

\bibitem{FK06}
S.G. Foss, D. Korshunov (2006).
Heavy tails in a multi-server queue.
\textit{Queueing Systems} \textbf{52 (1)},
31--48.

\bibitem{FK12}
S.G. Foss, D. Korshunov (2012).
On large delays in multi-server queues with heavy tails.
\textit{Math.\ Oper.\ Res.} \textbf{37 (2)},
201--218.

\bibitem{FKZ13}
S.G. Foss, D. Korshunov, S. Zachary (2013).
\textit{An Introduction to Heavy-Tailed and Subexponential Distributions.}
Springer.

\bibitem{FM14}
S.G. Foss, M. Miyazawa (2014).
Two-node fluid network with a heavy-tailed random input: the strong stability case.
\textit{J. Appl.\ Prob.} \textbf{51 (A)},
249--265.

\bibitem{GK96}
B.V. Gnedenko, V.Y. Korolev (1996).
\textit{Random Summation.}
CRC Press, Boca Raton, FL, USA.

\bibitem{JL99}
P.R. Jelenkovi\'c, A.A. Lazar (1999).
Asymptotic results for multiplexing subexponential on-off processes.
\textit{Adv.\ Appl.\ Prob.} \textbf{31},
394--421.

\bibitem{Kingman61}
J.F.C. Kingman (1961).
The single server queue in heavy traffic.
\textit{Math.\ Proc.\ Cambridge Phil.\ Soc.} \textbf{57 (4)},
902--904.

\bibitem{Kingman62}
J.F.C. Kingman (1962).
On queues in heavy traffic.
\textit{J. Roy.\ Stat.\ Soc.\ Series B (Methodological)} \textbf{24 (2)},
383--392.

\bibitem{Kingman65}
J.F.C. Kingman (1965).
The heavy traffic approximation in the theory of queues.
In: \textit{Proc.\ Symposium on Congestion Theory},
eds.\ W.L. Smith, W.E. Wilkinson (University of North Carolina Press, Chapel-Hill),
137--159.

\bibitem{LTWW94}
W.E. Leland, M.S. Taqqu, W. Willinger, D.V. Wilson (1994).
On the self-similar nature of Ethernet traffic (extended version).
\textit{IEEE/ACM Trans.\ Netw.} \textbf{2 (1)},
1--15.

\bibitem{LM08}
P.M.D. Lieshout, M.R.H. Mandjes (2008).
Asymptotic analysis of L\'evy-driven tandem queues.
\textit{Queueing Systems} \textbf{60 (3--4)},
203--226.

\bibitem{MB00}
M.R.H. Mandjes, S.C. Borst (2000).
Overflow behavior in queues with many long-tailed inputs.
\textit{Adv.\ Appl.\ Prob.} \textbf{32 (4)},
1150--1167.

\bibitem{NWZ22}
J.K. Nair, A. Wierman, B. Zwart (2022).
\textit{The Fundamentals of Heavy Tails: Properties, Emergence, and Estimation.}
Cambridge University Press.

\bibitem{Pakes75}
A.G. Pakes (1975).
On the tails of waiting-time distributions.
\textit{J. Appl.\ Prob.} \textbf{12},
555--564.

%\bibitem{PW00a}
%K. Park, W. Willinger (2000).
%Self-similar network traffic: an overview.
%In: K. Park, W. Willinger (eds.).
%\textit{Self-Similar Network Traffic and Performance Evaluation}.
%Wiley, New York,
%1--38.

%\bibitem{PW00}
%K. Park, W. Willinger (eds.) (2000).
%\textit{Self-Similar Network Traffic and Performance Evaluation},
%Wiley, New York.

\bibitem{PF95}
V. Paxson, S. Floyd (1995).
Wide area traffic: the failure of Poisson modeling.
\textit{IEEE/ACM Trans. Netw.} \textbf{3 (3)},
226--244.

\bibitem{RS00}
S. Resnick, G. Samorodnitsky (2000).
A heavy traffic approximation for workload processes with heavy tailed service requirements.
\textit{Mgmt.\ Sc.} \textbf{46},
1236--1248.

\bibitem{Veraverbeke77}
N. Veraverbeke (1977).
Asymptotic behaviour of Wiener-Hopf factors of a random walk.
\textit{Stoch.\ Proc.\ Appl.} \textbf{5},
27–37.

%\bibitem{Whitt74}
%W. Whitt (1974).
%Heavy traffic limit theorems for queues: a survey.
%In: \textit{Mathematical Methods in Queueing Theory},
%ed.\ A.B. Clarke (Springer, Berlin),
%307--350.

\bibitem{Whitt00}
W. Whitt (2000).
The impact of a heavy-tailed service-time distribution upon the M/GI/s waiting-time distribution.
\textit{Queueing Systems} \textbf{36},
71--87.

\bibitem{Whitt02}
W. Whitt (2002).
\textit{Stochastic-Process Limits: An Introduction to Stochastic-Process Limits and their Application to Queues.}
Springer.

\bibitem{Williams16}
R.J. Williams (2016).
Stochastic processing networks.
\textit{Ann.\ Rev.\ Stat.\ Appl.} \textbf{3},
323--345.

\bibitem{WTLW95}
W. Willinger, M.S. Taqqu, W.E. Leland, D.V. Wilson (1995).
Self-similarity in high-speed packet traffic: analysis and modeling of Ethernet traffic measurements.
\textit{Statistical Science} \textbf{10},
67--85.

\bibitem{WTSW97}
W. Willinger, M.S. Taqqu, R. Sherman, D.V. Wilson (1995).
Self-similarity through high-variability: statistical analysis of Ethernet LAN traffic at the source
level.
\textit{IEEE/ACM Trans.\ Netw.} \textbf{5 (1)},
71--86.

\bibitem{Zwart01}
B. Zwart (2001).
\textit{Queueing Systems with Heavy Tails}.
PhD Thesis, Eindhoven University of Technology.

%\bibitem{Zwart01}
%B. Zwart (2001).
%Tail asymptotics for the busy period in the GI/G/1 queue.
%\textit{Math.\ Oper.\ Res.} \textbf{26 (3)},
%485--493.

\bibitem{ZBM04}
B. Zwart, S.C. Borst, M. Mandjes (2004).
Exact asymptotics for fluid queues fed by multiple heavy-tailed on-off flows.
\textit{Ann.\ Appl.\ Prob.} \textbf{14},
903--957.

\end{thebibliography}
\end{document}